        %%%%%%%%%%%%%%%%%%%%%%%%%%%%%%%%%%%%%%%%%%%%%%%%%
        %%%                                           %%%
        %%      ON THE RADICAL OF BRAUER ALGEBRAS      %%
        %%                                             %%
        %%%            by Fabio Gavarini              %%%
        %%%%%                                       %%%%%
        %%%%%%%%%%%%%%%%%%%%%%%%%%%%%%%%%%%%%%%%%%%%%%%%%

\input amstex
\documentstyle{amsppt}
\input epsf

\magnification=\magstep1
\hsize=6.5truein
\vsize=9truein
%\hcorrection{0.9cm}
%\vcorrection{1cm}
\baselineskip=.15truein

\font \smallrm=cmr10 at 10truept
\font \smalllrm=cmr10 at 9truept

%
 
%

% "dangerous bend" sign
\def \llongrightarrow {\,\relbar\joinrel\relbar\joinrel\relbar\joinrel
\rightarrow\,}
\def \lllongrightarrow
{\,\relbar\joinrel\relbar\joinrel\relbar\joinrel\relbar\joinrel\relbar
\joinrel\relbar\joinrel\rightarrow\,}

\def \longtwoheadrightarrow {\,\relbar\joinrel\twoheadrightarrow\,}
\def \llongtwoheadrightarrow
{\,\relbar\joinrel\relbar\joinrel\relbar\joinrel\twoheadrightarrow\,}

\def \saltino {\vskip7pt}
\def \salto {\vskip9pt}
\def \ssalto {\vskip13pt}
\def \Salto {\vskip1.3cm}

\def \N {{\Bbb N}}
\def \Z {{\Bbb Z}}

\def \V {{\Bbb V}}

\def \Bfx {{\Cal B}_f^{(x)}}
\def \Bfor {{\Cal B}_f^{(n)}}
\def \Bfze {{\Cal B}_f^{(0)}}
\def \Bfuno {{\Cal B}_f^{(1)}}
\def \Bfsp {{\Cal B}_f^{(-2n)}}
\def \Sfx {{\Cal S}_f^{(x)}}
\def \Sfor {{\Cal S}_f^{(n)}}
\def \Sfsp {{\Cal S}_f^{(-2n)}}
\def \Ker {\text{\it Ker}}
\def \Rad {\text{\it Rad}}
\def \tas {\text{\tt tas}}
\def \bas {\text{\tt bas}}
\def \as {\text{\tt as}}
\def \d {\text{\rm \bf d}}
\def \h {\text{\rm \bf h}}
\def \endor {{End}_{O(V)} \! \left( V^{\otimes f} \right)}
\def \endsp {{End}_{Sp(W)} \! \left( W^{\otimes f} \right)}

\document

\topmatter

{\ } 

\vskip-33pt  

\hfill   {To appear in  {\sl Mathematische Zeitschrift}   
%   
% \/}  {\bf 52},  no.~3 (2002), 809--834
%   
 } 
\hskip19pt   {\ }  

\vskip41pt

\title
    ON THE RADICAL OF BRAUER ALGEBRAS
\endtitle

\author
   Fabio Gavarini
\endauthor

\affil
  Universit\`a degli Studi di Roma ``Tor Vergata'' ---
Dipartimento di Matematica  \\
  Via della Ricerca Scientifica 1, I-00133 Roma --- ITALY  
\endaffil

\address\hskip-\parindent
   Fabio Gavarini  \newline
   \indent   Universit\`a degli Studi di Roma ``Tor Vergata''
---   Dipartimento di Matematica  \newline
   \indent   Via della Ricerca Scientifica 1, I-00133 Roma, ITALY
---   e-mail: gavarini\@{}mat.uniroma2.it  \newline
   \indent   {\tt http://www.mat.uniroma2.it/\~{}gavarini/page-web.html}
\endaddress

%   \footnote""{ Revised version of December  $ 13^{\text{th}} $  2007 }
%
  \footnote""{ 2000 {\it Mathematics Subject Classification,}
Primary 20G05, 16N99; Secondary 16G99, 15A72. }
  \footnote""{ Partially supported by the European RTN ``LieGrits'',
contract no.~MRTN-CT-2003-505078, and by the Italian PRIN 2005
``Moduli e teorie di Lie''. }

 \abstract
   The radical of the Brauer algebra  $ \Bfx $  is known to be non-trivial
when the parameter  $ x $  is an integer subject to certain conditions (with
respect to  $ f \, $).  In these cases, we display a wide family of elements
in the radical, which are explicitly described by means of the diagrams of
the usual basis of  $ \Bfx $.  The proof is by direct approach for  $ \,
x = 0 \, $,  and  via classical Invariant Theory in the other cases,
exploiting then the well-known representation of Brauer algebras as
centralizer algebras of orthogonal or symplectic groups acting on
tensor powers of their standard representation.  This also gives a
great part of the radical of the generic indecomposable  $ \Bfx $--modules.
We conjecture that this part is indeed the whole radical in the case of
modules, and it is the whole part in a suitable step of the standard
filtration in the case of the algebra.  As an application, we find some
more precise results for the module of pointed chord diagrams, and for
the Temperley-Lieb algebra
   \hbox{--- realised inside  $ \Bfuno $  ---   acting on it.}
 \endabstract

\endtopmatter

 \vskip-21pt
\hfill  \hbox{\vbox{
     \hbox{\it \hskip22pt  ``Ahi quanto a dir che sia \`e cosa dura }
     \hbox{\it \hskip39pt     lo radical dell'algebra di Brauer }
     \hbox{\it \hskip19pt   pur se'l pensier gi\`a muove a congettura'' }
        \vskip4pt
     \hbox{\sl \hskip57pt     N.~Barbecue, ``Scholia'' } }
\hskip0,4truecm }

\Salto
 \vskip-21pt

\centerline { \bf  Introduction }

\salto

  The Brauer algebras first arose in Invariant Theory (cf.~[Br]) in
connection with the study of invariants of the action of the orthogonal
or the symplectic group   --- call it  $ G(U) $  ---   on the tensor powers
of its standard representation  $ U $.  More precisely, the centralizer
algebra  $ \, {End}_{G(U)} \big( U^{\otimes f} \big) \, $  of such an
action can be described by generators and relations: the latter depend
on the relationship among two integral parameters,  $ f $  and  $ x $   ---
the latter being related to  $ dim(U) $  ---   but when  $ x $  is big
enough (what is called ``the stable case'') the relations always remain
the same.  These ``stable'' relations then define an algebra  $ \Bfx $
of which the centralizer one is a quotient, obtained by adding the
further relations, when necessary. The abstract algebra  $ \Bfx $
is the one which bears the name of ``Brauer algebra''.
                                                \par
  The definition of  $ \Bfx $  still makes sense with  $ x $  arbitrarily
chosen in a fixed ground ring.  An alternative description is possible
too, by displaying an explicit basis of  $ \Bfx $  and assigning the
multiplication rules for elements in this basis.
                                                \par
  Assume the ground field  $ \Bbbk $  has characteristic zero.  Then
$ G(U) $  is linearly reductive, so by Schur duality the algebra  $ \,
{End}_{G(U)} \big( U^{\otimes f} \big) \, $  is semisimple: hence in the
stable case, when  $ \, \Bfx \cong {End}_{G(U)} \big( U^{\otimes f} \big)
\, $,  the Brauer algebra is semisimple too.  Otherwise,  $ \Bfx $  may
fail to be semisimple, i.e.~it may have a non-trivial radical.
                                                \par
  The most general result on  $ \, \Rad\left( \Bfx \right) \, $,  \, for
$ \, \text{\sl Char}\,(\Bbbk) = 0 \, $,  \, was found in [Wz]: for ``general
values'' of  $ x $   --- i.e., all those out of a finite range (depending on
$ f \, $,  \, and yielding the stable case) of values in the prime subring of
$ \Bbbk $  ---   the Brauer algebra  $ \Bfx $  is semisimple.  So the problem
only remained of computing  $ \Rad\left(\Bfx\right) $  when  $ x $  is an
integer and we are not in the stable case.  In this framework, the first
contributions came from Brown, who reduced the task to studying the radical
of ``generalized matrix algebras'' (cf.~[Br1--2]).  In particular, this
radical is strictly related with the nullspace of the matrix of structure
constants of such an algebra: later authors mainly followed the same strategy,
see e.g.~[HW1--2].  Further results were obtained using new techniques: see
[GL], [DHW], [KX], [CDM], [Hu], [DH].
                                                \par
  In the present paper we rather come back to the Invariant Theory
viewpoint.  The idea we start from is a very na\"\i{}ve one: as the
algebra  $ \, {End}_{G(U)} \big( U^{\otimes f} \big) \, $  is
semisimple, we have
 \vskip-8pt
  $$  \Rad \left( \Bfx \right) \;\; \subseteq \;\; \Ker\, \left( \pi_U :
\Bfx \! \longtwoheadrightarrow \! {End}_{G(U)} \big( U^{\otimes f} \big)
\right)  $$
 \vskip-8pt
\noindent
 where  $ \, \pi_U \! : \Bfx \!\!
\relbar\joinrel\relbar\joinrel\twoheadrightarrow \! {End}_{G(U)}
\big( U^{\otimes f} \big) \, $  is the natural epimorphism.  The second
step is an intermediate result,
%   
% of great interest in itself,   
%   
 namely
a description of the kernel  $ \, \Ker\, \left( \pi_U \! : \Bfx \!\!
\relbar\joinrel\relbar\joinrel\twoheadrightarrow \! {End}_{G(U)} \big(
U^{\otimes f} \big) \right) \, $.  Indeed, using the Second Fundamental
Theorem of classical invariant theory we find a set of linear generators
for it: they are explicitly written in terms of the basis of diagrams,
and called  {\it (diagrammatic) minors  {\rm or}  Pfaffians},  depending
on the sign of  $ x \, $.  As  $ \Ker\,\big(\pi_U\big) $  contains
$ \Rad \left( \Bfx \right) $,  every element of  $ \Rad\left(\Bfx\right) $
is a linear combination of these special elements (minors or Pfaffians).
As a last step, a basic knowledge of  $ \Bfx $--modules
yields some more information on the structure of the semisimple quotient
of  $ \Bfx $.  Thus we determine exactly which ones among minors, or
Pfaffians, belong to  $ \Rad\left( \Bfx \right) \, $:  \, so we find
a great part of  $ \Rad\left( \Bfx \right) \, $,  \, and we conjecture
that this is  {\sl all\/}  the part of the radical inside the proper
step of the standard filtration.  We then find a similar result and
conjecture for the generic indecomposable  $ \Bfx $--modules  too.
                                                \par
  Our approach applies directly only in case  $ x $  is an integer which
is not zero nor odd negative; but  {\sl a posteriori},  we find also
similar results for  $ \, x = 0 \, $,  via an  {\sl ad hoc\/}  approach.
                                                \par
  Also, we discuss how much of these results can be extended to the case
of  $ \text{\sl Char}\,(\Bbbk) > 0 \, $.
 \vskip3pt
  Finally, we provide some more precise results for the module of pointed
chord diagrams, and the Temperley-Lieb algebra   --- realised as a
subalgebra of  $ \Bfuno $  ---   acting on it.

\vskip5pt

\centerline{ \smalllrm ACKNOWLEDGEMENTS }
 \vskip-1pt
   \centerline{\smallrm The author thanks A.~Knutson, O.~Mathieu,
G.~Papadopoulo, P.~Papi and C.~Procesi}
 \vskip-2pt
   \centerline{\smallrm for several useful discussions,
and Jun Hu and Hebing Rui for their valuable remarks.}
 \Salto

\centerline{ \bf  \S 1  The Brauer algebra }

\ssalto

   {\bf 1.1  $ f $--diagrams.} \  Given  $ \, f \in \N_+ \, $,  \, denote
by  $ \V_f $  the datum of  $ 2f $  vertices in a plane, arranged in
two rows, one upon the other, each one of  $ f $  aligned vertices.  Then
consider the graphs with  $ \V_f $  as set of vertices and  $ f $  edges,
such that each vertex belongs to exactly one edge.
%
% The picture below
% shows an example of such a graph for  $ \, f=6 \, $.
% %
% %   PICTURE OF A 6-DIAGRAM
% %
% \vskip3pt
%   $$  \epsfbox{6-diagram.eps}  $$
% \vskip3pt
% %
% %
%
   We call such graphs  {\it  $ f $--diagrams\/},  denoting by  $ D_f \, $
the set of all of them.  In general, we shall denote them by bold roman
letters, like  $ \d $.  These  $ f $--diagrams are as many as the pairings
of  $ 2f $  elements, hence  $ \, \big| D_f \big| = (2f-1)!! := (2f-1)
\cdot (2f-3) \cdots 5 \cdot 3 \cdot 1 \, $  in number.
                                               \par
   We shall label the vertices in  $ \V_f $  in two ways: either we label
the vertices in the top row with the numbers  $ 1^+ $,  $ 2^+ $,  $ \dots $,
$ f^+ $,  in their natural order from left to right, and the vertices in the
bottom row with the numbers  $ 1^- $,  $ 2^- $,  $ \dots $,  $ f^- $,  again
from left to right, or we label them by setting  $ i $  for  $ i^+ $  and
$ f+j $  for  $ j^- $  (for all  $ \, i, j \in \{1,2,\dots, f\} $).  Thus
an  $ f $--diagram  can also be described by specifying its set of edges.
%
% :
% for instance, the  $ 6 $--diagram  above is given by
% $ \, \big\{ \{1^+,4^+\}, \{3^-,5^+\}, \{2^+,4^-\}, \{5^-,6^+\}, \{2^-,6^-\},
% \{3^+,1^-\} \big\} \, $.
%
 In general, given  $ f $--tuples  $ \, \text{\rm
\bf i} := (i_1, i_2, \dots, i_f) $  and  $ \text{\rm \bf j} := (j_1, j_2,
\dots, j_f) $  such that  $ \, \{i_1, \dots, i_f\} \cup \{j_1, \dots,
j_f\} = \V_f \, $,  we call  $ \, \d_{\text{{\bf i}, {\bf j}}} \, $  the
$ f $--diagram  obtained by joining  $ i_k $  to  $ j_k $, for each  $ \,
k= 1, 2, \dots, f \, $.
%
% For instance, the above diagram is  $ \d_{\text{
% {\bf i}, {\bf j}}} $  for  $ \, \text{\bf i} = \{1^+, 2^+, 3^+, 5^+, 6^+,
% 2^-\} \, $,  $ \, \text{\bf j} = \{4^+, 4^-, 1^-, 3^-, 5^-, 6^-\} \, $.
%
                                                        \par
   When looking at the edges of an  $ f $--diagram,  we shall distinguish
between those which link two vertices in the same (top or bottom) row,
which we call  {\it horizontal edges}  or  simply  {\it arcs},  and those
which link two vertices in different rows, to be called  {\it vertical
edges}.  Clearly, any  $ f $--diagram  $ \d $  has the same number of arcs
in the top row and in the bottom row: if this number is  $ k $,  we shall
say that  $ \d $  is a  $ k $--arc  ($ f $--)diagram.  Then, letting
$ \; D_{f,k} := \big\{\, \d \in D_f \,\big\vert\, \d \hbox{\it \ is
a  $ k $--arc  diagram} \,\big\} \; $  we have  $ \; D_f =
\bigcup_{k=1}^{[f/2]} D_{f,k} \; $;  \; hereafter, for any
$ \, f \in \N \, $  we set  $ \, \big[f/2\big] := f/2 \; $  if
$ f $  is even and  $ \, \big[f/2\big] := (f\!-\!1)\big/2 \; $
if  $ f $  is odd.
%
%  such notation as  $ \big[f/2\big] $  denotes the  {\sl integral
% part\/}  (also called the  {\sl ``floor''\/})  of  $ f/2 \, $,
% for any  $ \, f \in \N \, $.
%

\salto

   {\bf 1.2  Arc structure and permutation structure of diagrams.} \  Let
$ \d $  be an  $ f $--diagram.  With ``arc structure of the top row'',
resp.~``bottom row'', of  $ \d $  we shall mean the datum of the arcs in
the top, resp.~bottom, row of  $ \d $,  in their mutual positions.  To put
it in a nutshell, we shall use such terminology as ``top arc structure'',
resp.~``bottom arc structure'',  of  $ \d $   --- to be denoted with
$ \, \tas(\d) \, $,  resp.~$ \, \bas(\d) $  ---   and ``arc structure of
$ \d $''   --- to be denoted with  $ \, \as(\d) $  ---   to mean the datum
of both top and bottom arc structures of  $ \d $,  that is  $ \, \as(\d)
:= \big( \tas(\d), \bas(\d) \big) \, $.  Note that any top or bottom arc
structure may be described by a one-row graph of vertices, arranged on a
horizontal line, and some edges (the arcs) joining them pairwise, so that
each vertex belongs to at most one edge.  Following Kerov (cf.~[Ke]), such
a graph will be called a  {\it  $ k $--arc  $ f $--junction},  or  {\it
$ (f,k) $--junction},  where  $ f $  is its number of vertices and  $ k $
its number of edges.
                                                      \par
   We denote the set of  $ (f,k) $--junctions  by  $ \, J_{f,k}
\, $.  Then clearly  $ \, \big\vert J_{f,k} \big\vert = \Big(\!
{f \atop 2k} \Big) (2k-1)!! \, $.

\saltino

  Any  $ \, \d \in D_{f,k} \, $  has exactly  $ f-2k $  vertices in its top
row, and  $ f-2k $  vertices in its bottom row which are pairwise joined
by its  $ f-2k $ vertical edges.  Let us label with  $ 1 $,  $ 2 $,
$ \dots $,  $ f-2k $  from left to right the vertices in the top row,
and do the same in the bottom row.  Then there exists a unique permutation
$ \, \sigma = \sigma(\d) \in S_{f-2k} \, $   --- to be called the
``permutation structure'', or ``symmetric (group) part'', of  $ \d $
---   such that  $ \sigma(i) $  is the label of the bottom row vertex
of  the vertical edge whose top row vertex is labelled with  $ i \, $.
                                               \par
   Therefore the maps  $ \, \d \mapsto \big(\sigma(\d), \as(\d) \big)
\, $  set bijections  $ \; D_{f,k}  \llongrightarrow  S_{f-2k} \times
\big( J_{f,k} \times J_{f,k} \big) \; $
%
%   \eqno (1.1)  $$
%
 and altogether they give a bijection
$ \, D_f \llongrightarrow \bigcup_{k=1}^{[f/2]} S_{f-2k} \times
\left( {J_{f,k}}^{\!\times 2} \right) \, $.

\salto

   {\bf 1.3 The Brauer algebra.} \  Let  $ \Bbbk $  be a field,  $ \, p :=
\text{\sl Char}\,(\Bbbk) \geq 0 \, $,  \, and take  $ \, x \in \Bbbk \, $.
Later results will require some restrictions on  $ \Bbbk \, $,  \, but on
the other hand one can also generalise, replacing  $ \Bbbk $  with any
commutative ring with 1   --- see (4.1) and the subsequent remark.
                                                         \par
   Let  $ \Bfx $  be the  $ \Bbbk $--vector  space with basis  $ D_f
\, $;  we introduce a product in  $ \Bfx $  by defining the product
of  $ f $--diagrams  and extending by linearity.  So for all  $ \,
\text{\bf a}, \text{\bf b} \in D_f \, $  define the product  $ \,
\text{\bf a} \cdot \text{\bf b} = \text{\bf a} \text{\bf b} \, $  as
follows.  First, draw  $ \text{\bf b} $  below  $ \text{\bf a} $;  second,
connect the  $ i $--th  bottom vertex of  $ \text{\bf a} $  with the
$ i $--th  top vertex of  $ \text{\bf b} \, $;  third, let  $ C(\text{\bf a},
\text{\bf b}) $  be the number of cycles in the new graph obtained in (2)
and let  $ \text{\bf a} \ast \text{\bf b} $  be this graph,  {\sl pruning
out the cycles\/};  then  $ \text{\bf a} \ast \text{\bf b} $  is a new
$ f $--diagram,  and we set  $ \; \text{\bf a} \, \text{\bf b} :=
x^{C(\text{\bf a}, \text{\bf b})} \text{\bf a} \ast \text{\bf b} \; $.
%
% We denote by  $ \, \ast : D_f \times D_f \rightarrow D_f \, $  the map
% given by  $ \, (\text{\bf a}, \text{\bf b}) \mapsto \text{\bf a} \ast
% \text{\bf b} \, $  and  $ \, C : D_f \times D_f \rightarrow \N \, $  the
% map given by  $ \, (\text{\bf a}, \text{\bf b}) \mapsto C(\text{\bf a},
% \text{\bf b}) \, $.
%
% The following is an example:
% %
% %     PICTURE OF THE PRODUCT OF TWO DIAGRAMS
% %
% \vskip5pt
%   $$  \epsfbox{prod_ab-c.eps}  $$
% \vskip-86,5pt
%   $$  {} \hskip30pt  \phantom{\hbox{\bigrm = \hskip13pt \bigit x}}
% \hskip28pt {}  $$
% \vskip59pt
%
%
  This definition makes  $ \Bfx $  into a
a unital associative  $ \Bbbk $--algebra:  this is the
{\sl Brauer algebra},  in its ``abstract'' form (see for instance
[KX]).  Its relation with Brauer's centralizer algebra is explained
in \S 4 later on.
                                                         \par
  Note that for  $ \text{\bf a}, \text{\bf b} \in D_f \, $,  the
top, resp.~bottom, arc structure of  $ \text{\bf a} \ast \text{\bf b} $
``contains'' that of  $ \text{\bf a} $,  resp.~$ \text{\bf  b} \, $.  In
particular, if  $ \, \text{\bf a} \in D_{f,a} \, $  and  $ \, \text{\bf b}
\in D_{f,b} \, $  this gives  $ \, \text{\bf a} \ast \text{\bf b} \in
D_{f,\max(a,b)} \, $.
                                                         \par
  The symmetric group  $ S_{2f} $  acts on  $ \V_{\!f} \, $,  once a numbering
of the vertices in  $ \V_{\!f} $  is fixed; so it acts on  $ D_f \, $,  and
linear extension gives a  $ \Bbbk $--linear  action on  $ \Bfx $  (studied
in [DH], [Hu]).
                                                         \par
   By construction  $ D_{f,0} $  is a subset of  $ \Bfx $,  actually a
subsemigroup.  Now, for any  $ \, \sigma \in S_f \, $  let  $ \, \d_\sigma
\in D_{f,0} \, $  be the  $ f $--diagram  obtained by joining  $ i^+ $  with
$ {\sigma(i)}^- $  (cf.~\S 1.1).  Then the map  $ \, S_f \rightarrow D_{f,0}
\subset \Bfx \, $  is a morphism of semigroups, whose image is  $ D_{f,0} \, $;
\, hence  $ \Bfx $  contains a copy of  $ S_f $  (namely  $ D_{f,0} $)  and a
copy of the group algebra  $ \Bbbk \left[ S_f \right] \, $.  Restricting the
left (right) regular representation of  $ \Bfx $  we get a left (right) action
of  $ S_f $  on  $ \Bfx \, $.

\salto

   {\bf 1.4  Presentation of  $ \Bfx $,  and signs of diagrams.} \  By
\S 1.3,  $ \Bfx $  contains a copy of the symmetric group  $ S_f \, $.
Moreover, for any pair of distinct indices  $ \, i, j \in \{1,2,\dots,f\}
\, $  we define  $ \, \h_{i,j} \in D_{f,1} \, $  to be the  $ f $--diagram
with an arc joining  $ i^+ $  with  $ j^+ $,  an arc joining  $ i^- $  with
$ j^- $,  and a vertical edge joining  $ k^+ $  with  $ k^- $  for every
$ \, k \in \{1,2,\dots,f\} \setminus \{i,j\} \, $.
                                                   \par
   The  $ \h_{i,j} $'s  together with the elements  $ \d_\sigma \in D_{f,0}
\, $  (for all  $ \, \sigma \in S_f \, $)  generate the algebra  $ \Bfx \, $;
\, the relations among these generators are known too (see, e.g., [DP], \S 7).
As  $ D_{f,1} $  is a single  $ D_{f,0} $--orbit  (i.e.~$ S_f $--orbit),
taking only one  $ 1 $--arc  $ f $--diagram  is enough, so  $ \Bfx $  is
generated, for instance, by  $ D_{f,0} \bigcup \, \{\, \h_{1,2} \} \, $.
In particular, for any  $ \, \d \in D_{f,k} \, $  there exist  {\sl unique\/}
$ \, \d_\sigma, \d_\rho \in D_{f,0} \, $  such that  $ \, \d = \d_\sigma \,
\h_{1,2} \cdots \h_{2k-1,2k} \, \d_\rho \, $  and moreover  $ \sigma $  and
$ \rho $  do not invert any of the pairs  $ (1,2) $,  $ (3,4) $,  $ \dots $,
$ (2k-1,2k\,) \, $.  Then we define the  {\it sign of\/}  $ \d $  to be  $ \,
\varepsilon(\d):= sgn(\sigma) \cdot {(-1)}^k \cdot sgn(\rho) \, $,  \, which
is independent of the given factorization of  $ \d \, $.

\Salto
\vskip-5pt

\centerline{ \bf  \S 2  The standard series and  $ \Bfx $--modules }

\ssalto

   {\bf 2.1  The standard series.} \  For any  $ \, k \in \{1, 2, \dots
[f/2] \} \, $,  let  $ \, \Bfx \langle k \rangle \, $  be the vector
subspace of  $ \Bfx $  spanned by  $ D_{f,k} \, $.  We define  $ \,
\Bfx(k) := \bigoplus_{h=k}^{[f/2]} \Bfx \langle h \rangle \, $,
\, so  $ \Bfx(k) $  has  $ \Bbbk $--basis  $ \bigcup_{h=k}^{[f/2]}
D_{f,h} \; $.  By definition, the  $ \Bfx(k) $'s  form a chain of
subspaces
  $$  \Bfx = \, \Bfx(0) \, \supset \, \Bfx(1) \, \supset \, \cdots \,
\supset \, \Bfx(k) \, \supset \, \cdots \, \supset \, \Bfx \big( [f/2]
\big) \, \supset \, \big\{0\big\}   \eqno (2.1)  $$
 \eject
\noindent
 which we call ``standard series''.  We denote by  $ \,  \Bfx[k] := \Bfx(k)
\Big/ \Bfx(k+1) \, $  the  $ (k + \! 1) $-th   factor (a quotient space)
of this series, setting also  $ \; \Bfx \big( [f/2]+1 \big) \, :=
\big\{0\big\} \; $.

\saltino

   By construction each  $ \Bfx(k) $  is a (two-sided) ideal of
$ \Bfx $:  thus every  $ \Bfx[k] $  inherits a structure of associative
$ \Bbbk $--algebra,  one of left  $ \Bfx $--module,  and one of right
$ \Bfx $--module.  Moreover  $ \, \Bfx\!(k) = \Bfx \! \langle k \rangle \,
\oplus \, \Bfx\!(k+1) \, $,  so any basis for  $ \Bfx \! \langle k \rangle $,
taken modulo  $ \Bfx\!(k+1) $,  serves as basis for  $ \Bfx[k] \, $;  in
particular, we shall use  $ D_{f,k} $  as a basis of $ \Bfx[k] \, $.  Note
that, since the  $ \Bfx(k) $'s  are two-sided ideals of  $ \Bfx $,  the
$ \Bfx[k] $'s  are  $ \Bfx $--bimodules  too.
%
%                                                              \par
%   The  $ \Bfx(k) $'s  and the  $ \Bfx[k] $'s  have been studied extensively
% in [Bw2]   --- where the algebra structure of the  $ \Bfx[k] $  is
% determined ---   and in [HW1--2]   --- where the link between the
% representations of  $ \Bfx $  and the structure of the  $ \Bfx[k] $'s
% is explained, and the study of the radical of  $ \Bfx $  is reduced to
% the study of the radicals of the  $ \Bfx[k] $'s.  More recently, it was
% proved in [KX] that each  $ \Bfx[k] $  is a so-called  {\sl inflation
% algebra},  so that every  $ \Bfx(k) $   --- hence also  $ \Bfx $  itself
% ---   is an  {\sl iterated inflation algebra}.
%                                                               \par
%    In the next section we recall some of Brown's results.
%

\salto

   {\bf 2.2  The structure of  $ \Bfx[k] \, $.} \  Let us fix some more
notation.  Given  $ \, h \in \N \, $,  we write  $ \, \lambda \vdash
h \, $  to mean that  $ \lambda $  is a partition of  $ h \, $;  then
for  $ \, \lambda \vdash h \, $  we denote by  $ \lambda^t $  the dual
partition.  Also, if  $ \, \lambda \vdash h \, $  we denote by
$ M_\lambda $  the unique associated simple  $ S_h $--module,  with
the assumption that  $ M_{(h)} $  is the trivial  $ S_h $--module
and  $ M_{(\undersetbrace{h}\to{\scriptstyle 1,1,\dots,1})} $  is
the alternating one.
                                                    \par
   Now let  $ \, k \in \{1, 2, \dots [f/2] \} \, $  be fixed.  Consider
the set  $ \, J_{f,k} \, $  of  $ (f,k) $--junctions  defined in \S 1.2,
and define  $ \, H_{f,k} \, $  to be the  $ \Bbbk $--vector  space with
basis  $ J_{f,k} \; $.  In particular, one has that  $ \, dim(H_{f,k})
= \big\vert J_{f,k} \big\vert = \Big(\! {f \atop 2k} \!\Big) (2k-1)!!
\; $.  Inverting the map  $ \, D_{f,k} \longrightarrow S_{f-2k} \times
\big( J_{f,k} \times J_{f,k} \big) \, $ (cf.~\S 1.2) and extending by
linearity we define two linear isomorphisms
 \vskip-5pt   
  $$  \eqalign{
   \langle \boxtimes \rangle  &  : \, \Bbbk \big[ S_{f-2k} \big]
\otimes \big( H_{f,k} \otimes H_{f,k} \big)  \llongrightarrow
\Bfx \langle k \rangle  \cr
   \boxtimes  &  : \, \Bbbk \big[ S_{f-2k} \big] \otimes
\big( H_{f,k} \otimes H_{f,k} \big)  \llongrightarrow  \Bfx[k]  \cr }
\eqno (2.2)  $$
 \vskip-1pt   
   By Young's theory,  $ \Bbbk \big[ S_{f-2k} \big] $  splits into  $ \;
\Bbbk \big[ S_{f-2k} \big] = \! \bigoplus\limits_{\mu \vdash (f-2k)} \!
I_\mu \; $.  Hereafter, each  $ I_\mu $  is a two-sided ideal of  $ \Bbbk
\big[ S_f \big] $,  and a simple algebra, namely the algebra of linear
endomorphisms of the simple  $ S_{f-2k} $--module  $ M_\mu $,  which is
a full matrix algebra over  $ \Bbbk \, $.  Then we set

\salto

\proclaim{Definition 2.3}  For every  $ \, \mu \vdash (f-2k) \, $  we
define  $ \; \Bfx[k \, ; \mu] := \boxtimes \big( I_\mu \otimes (H_{f,k}
\otimes H_{f,k}) \big) \; $  and  $ \; \Bfx \langle k \, ; \mu \rangle :=
\langle \boxtimes \rangle \big( I_\mu \otimes (H_{f,k} \otimes H_{f,k})
\big) \; $.  Moreover, we denote by  $ \; \Bfx(k \, ; \mu) \; $  the
preimage of  $ \, \Bfx[k \, ; \mu] $  in  $ \, \Bfx(k) \, $.
\endproclaim

\salto

   {\bf 2.4 Generalized matrix algebras.} \  We recall (from [Bw1]) the notion of  {\sl generalized matrix algebra\/}:  this is any associative
$ \Bbbk $--algebra  $ A $  with a finite basis  $ \, \{\, e_{ij} \,\}_{i,j \in I} \, $  for which the multiplication table looks like  $ \; e_{ij} \cdot e_{pq} = \sigma^\ast_{jp} \, e_{iq} \; $  for some  $ \, \sigma^\ast_{jp} \in \Bbbk \;\, \big(\, \forall \; i, j, p, q \in I \,\big) \, $.  
% Then we denote by  $ \, \Phi(A) := {\big\{\, \sigma^\ast_{ij} \big\}}_{i,j \in I} \, $  the corresponding matrix of structure constants.    
%                                                            \par
%    If  $ A $  is any generalized matrix algebra, then the following hold (cf.~[Bw1]):
% 
Then we set  $ \, \Phi(A) := {\big\{\, \sigma^\ast_{ij} \big\}}_{i,j \in I} \, $.  For such an  $ A $,  the following hold (cf.~[Bw1]):
  \roster
     \item  either  $ A $  is simple, or  $ A $ has non-zero radical
$ \Rad\,(A) \, $,  \, and  $ A \big/ \Rad\,(A) $  is simple;
     \item  $ A $  is simple if and only if it has an identity element;
     \item  $ dim_\Bbbk(A) = h^2 \, $,  \, for some  $ \, h \in \N \, $,  \, and  $ \; dim_\Bbbk \big( \Rad\,(A) \big) = h^2 - {rk \big( \Phi(A)
\big)}^2 \; $;
     \item  the nilpotency degree of  $ \Rad\,(A) $  is at most  $ 3 \, $.
  \endroster

\vskip2pt   
  
\noindent   
 Hereafter,  {\sl by ``radical''  $ \Rad\,(\frak{A}) $  of any (possibly non-unital) algebra\/}  $ \frak{A} $  we shall mean the intersection of the annihilators of all its simple left modules (Brown's definition is the set of properly nilpotent elements: for generalized matrix algebras, the two definitions coincide).   

\vskip3pt   

   The most general result about the structure of  $ \Bfx[k] $  is the next one:

\salto

\proclaim{Theorem 2.5 {\sl ([Bw2], \S\S 2.2--3; [KX], \S\S 3--5)}}  For any
$ \, \mu \vdash (f\!-\!2k) \, $,  \, the subspace  $ \, \Bfx [k\,;\mu] $
is a two-sided ideal of  $ \, \Bfx[k] \; $,  \, the algebra  $ \, \Bfx[k] $
splits as
%
% a direct sum of these ideals,
%
 $ \;\; \Bfx[k] = \bigoplus_{\mu \vdash (f-2k)} \Bfx[k\,;\mu] $
%
% \eqno (2.3)  $$
%
%    \indent
%
\; and the  $ \, \Bfx[k\,;\mu] $'s
%
% (for  $ \, \mu \vdash (f\!-\!2k) \, $)
%
are pairwise non-isomorphic generalized matrix algebras.

 Moreover, every  $ \, \Bfx(k\,;\mu) \, $  is a two-sided ideal of  $ \;
\Bfx(k) \, $,
 \; and every  $ \, \Bfx[k\,;\mu] \, $  is a  $ \Bfx $--sub-bimodule
of  $ \, \Bfx[k] \, $,
\; for any  $ \, \mu \vdash (f\!-\!2k) \; $.
\endproclaim

%
% \salto
%
%   A further fact, easily deduced from [Bw2], is the following:
%
% \salto
%
% \proclaim{Lemma 2.6}  Every  $ \, \Bfx[k\,;\mu] \, $,  \; for  $ \, \mu
% \vdash (f-2k) \, $,  \, is a  $ \Bfx $--sub-bimodule  of  $ \, \Bfx[k]
% \; $.
% \endproclaim
%

\saltino

   {\bf 2.6  Representations of  $ \Bfx $.} \  Let  $ \, 0 \leq k \leq
f \big/ 2 \, $,  and let  $ H_{f,k} $  be the vector space defined in
\S 2.2.  For any  $ \, \mu \vdash (f-2k) \, $  ($ k \in \{0, 1, \dots,
[f/2]\} $)  we define  $ \, H^\mu_{f,k} := M_\mu \otimes H_{f,k} \; $.
We endow  $ H^\mu_{f,k} $  with a structure of  $ \Bfx $--module,
following Kerov (cf.~[Ke], [HW1--2], [GP]).
                                            \par
   Let  $ \d $  be an  $ f $--diagram,  and let  $ v $  be an
$ (f,k) $--junction.  For all  $ i=1,\ldots,f $,  connect the  $ i $--th
bottom vertex of  $ \d $  with the  $ i $--th  vertex of  $ v $,  let
$ C(\text{\bf d},v) $  be the number of loops occurring in the new graph
$ \Gamma(\d,v) $  obtained in this way, and let  $ \, a \star v \, $  be
the graph made of the vertices of the top line of  $ \d $,  connected
by an edge iff they are connected (by an edge or a path) in the new graph
$ \Gamma(\d,v) \, $.  Then  $ \, \d \star v \in J_{f,k'} \, $,  with
$ k' \geq k $  and  $ \, k' = k \, $  iff each pair of vertices of
$ v $  which are connected by a path in  $ \Gamma(\d,v) $  are in fact
joined by an edge in  $ v \, $:  in this case we say that the junction
$ v $  is {\it admissible}  for the diagram  $ \d $.  We set
 \vskip-7pt
  $$  \d.v \, := \, x^{C(\text{\bf d},v)} \, \d \star v  \quad
\hbox{\sl if  $ v $  is admissible for  $ \d $} \, ,   \qquad
\d.v \, := \, 0  \quad  \hbox{\sl otherwise} \, .  $$
%
%    \indent   The following are two, simple examples:
% %
% %   PICTURE OF \d.v = x^c \, \d \star v
% %
% \vskip7pt
%   $$  \epsfbox{diag_junct1.eps}  $$
% \vskip-68pt
%   $$  {} \hskip29pt  \phantom{\hbox{\bigrm = \hskip14pt \bigit x}}
% \hskip28pt {}  $$
% %
%  \vskip31pt
% %
% %
% %   PICTURE OF \d.v = 0
% %
%   $$  {} \hskip-77pt  \epsfbox{diag_junct2.eps}  $$
% %
%   $$  {} \hskip211pt  \phantom{\hbox{\bigrm = \hskip14pt 0}}
% \hskip28pt {}  $$
% %
% \vskip-11pt
% %
% %
%
 \vskip1pt
 See [Ga], \S 2.11 for some simple examples.
%
%    \indent
%                                           \par

 \saltino

 To any pair  $ \, (\d,v) \in D_f \times J_{f,k} \, $  we can
also attach an  element  $ \, \pi(\d,v) \in S_{f-2k} \, $:  this is the
permutation which carries   --- through the graph  $ \Gamma(\d,v) $  ---
the isolated vertices of  $ v $  into the isolated vertices of  $ \d \star
v $  (one takes into account only the relative position of the isolated
vertices in  $ v $,  $ \d \star v $)  in case  $ v $  is admissible for
$ a $,  otherwise it is  $ id \, $.
%
% In the previous examples we have  $ \,
% \pi(\d,v) = \left( {1 \ 2 \ 3 \atop 2 \ 1 \ 3} \right) \, $  in the
% non-trivial case.
%

\salto

\proclaim{Proposition 2.7  {\sl (cf.~[Ke], [HW1], [KX], [CDM])}}  Assume
$ \, \text{\sl Char}\,(\Bbbk) = 0 \, $  or  $ \, \text{\sl Char}\,(\Bbbk)
> f \, $.
                                   \par
   (a) \, Linear extension of the rule  $ \; \d.(u \otimes v):=
\pi(\d,v).u \,\otimes\, \d.v \; $  for every  $ \, (\d,v) \in D_f \times
J_{f,k} \, $  endows  $ H^\mu_{f,k} $  with a well-defined structure
of module over  $ \, \Bfx $.  Then  $ H^\mu_{f,k} $  is also a module
over  $ \, \Bfx \big/ \Bfx(k+1) \, $  and over  $ \, \Bfx[k] :=
\Bfx(k) \big/ \Bfx(k+1) \; $.
                                   \par
   (b) \, The various modules  $ H^\mu_{f,k} \, $,  for different pairs
$ (k,\mu) $   --- over any of the previous algebras ---   are pairwise
non-isomorphic.
                                   \par
   (c) \, If  $ \, \Bfx $  is semisimple, then every  $ \Bfx $--module
$ H^\mu_{f,k} $  is simple and, conversely, any simple  $ \Bfx $--module
is isomorphic to one of the  $ H^\mu_{f,k} $'s.
                                   \par
   (d) \, If every  $ \Bfx $--module  $ H^\mu_{f,k} $  is simple, then
the algebra  $ \, \Bfx $  is semisimple.
                                   \par
   (e) \, Every simple  $ \Bfx $--module  is a quotient of some
$ H^\mu_{f,k} \, $.  Conversely, each  $ H^\mu_{f,k} $  has a simple
quotient, but for the case of even  $ f \, $,  $ \, k = f \big/ 2 \, $
and  $ \, x = 0 \, $.  Indeed, these simple quotients are in bijection
with the isoclasses of simple  $ \, \Bfx $--modules.
\endproclaim

\Salto
\vskip-5pt

\centerline{ \bf  \S 3  Semisimple quotients of  $ \Bfx $
and of  $ \Bfx $--modules }

\ssalto

%
%   In this section we present some results about the semisimple quotients
% of the Brauer algebra and of its representations introduced in \S 2.
% Many of these facts are more or less known among specialists (see,
% for instance, [DHW] and references therein).  Nevertheless, we
% (re)formulate them for the sake of completeness and clarity, and
% --- more important ---   to fill in some gaps regarding the case
% $ \, x = 0 \, $.  We mostly base on Brown's work.
%
% \salto
%

   {\bf 3.1  Splitting the semisimple quotient of  $ \Bfx \, $.} \
Let  $ \; \Rad \left( \Bfx \right) \, $,  \; resp.~$ \; \Sfx := \,
\Bfx \! \bigg/ \! \Rad \left( \Bfx \right) \, $,  \; resp.~$ \; \pi_1
\, \colon \, \Bfx \!\relbar\joinrel\relbar\joinrel\twoheadrightarrow \Sfx
\; $  denote respectively the radical of  $ \Bfx $,  its semisimple quotient,
and the canonical epimorphism.  By general theory,  $ \Sfx $  has a direct
sum decomposition  $ \; \Sfx \, = \, \bigoplus_{i \in I} S_i \; $  in which
the  $ S_i $'s  are two-sided ideals which are simple algebras.  Of course,
$ \Sfx $  is a  $ \Bfx $--left/right/bi-module,  so the  $ S_i $'s  are
left/right/bi-submodules over  $ \Bfx $,  and each  $ S_i $  is simple as
a  $ \Bfx $-bimodule.  In this section we collect some information about
what the set  $ I $  has to be and what the blocks  $ S_i $'s  arise from.
                                              \par
   Define  $ \, \Sfx(k) := \pi_1 \left( \Bfx(k) \right) \, $  for all
$ \, k= 1, 2, \dots, [f/2] \, $:  then each  $ \Sfx(k) $  is a two-sided
ideal of  $ \Sfx $,  hence also a left/right/bi-submodule over  $ \Bfx
\, $.  In particular, there exists  $ I_k \subseteq I $  such that  $ \;
\Sfx(k) = \bigoplus_{i \in I_k} S_i \; $.  Applying  $ \pi_1 $  to (2.1),
one gets a series
  $$  \Sfx = \, \Sfx(0) \, \supseteq \, \Sfx(1) \, \supseteq \, \cdots \,
\supseteq \, \Sfx(k) \, \supseteq \, \cdots \, \supseteq \, \Sfx([f/2])
\, \supseteq \, \big\{0\big\}  $$
which corresponds to the chain of inclusions  $ \, I = I_0 \supseteq I_1 \supseteq \cdots \supseteq I_k \supseteq \cdots \supseteq I_{[f/2]} \supseteq \emptyset \; $.  So the algebra  $ \, \Sfx[k] := \Sfx(k) \Big/ \Sfx(k+1)
\, $  (for all  $ \, k \, $,  setting also  $ \, \Sfx \big( [f/2]+1 \big) := 0 \, $)  is well defined, and splits (up to isomorphisms) as  $ \, \Sfx[k] = \bigoplus_{j \in I_k \setminus I_{k+1}} \! S_j \; $,  \, a direct sum of simple algebras.  Moreover (up to isomorphisms),  $ \, \Sfx = \bigoplus_{k=0}^{[f/2]} \Sfx[k] \; $.
                                              \par
%
%    By construction, there
%
 There
is an algebra epimorphism  $ \, \pi_1^* \, \colon \,
\Bfx[k] \! \longtwoheadrightarrow \! \Sfx[k] \, $
%
% such that
%   $$  \CD
%    \Bfx(k)  @>{\; \pi_1 \;}>>  \Sfx(k)  \\
%   @V{p_B}VV                   @VV{p_S}V  \\
%    \Bfx[k]  @>>{\; \pi_1^* \;}>  \Sfx[k]  \\
%       \endCD  $$
% %
% (where  $ p_B $  and  $ p_S $  are the canonical projections) is a
%
 which together with  $ \pi_1 $  and the canonical projections from
$ \Bfx(k) $  to  $ \Bfx[k] $  and from  $ \Sfx(k) $  to  $ \Sfx[k] $
forms a commuta\-tive diagram.  Finally, define  $ \, \Sfx[k\,;\mu] :=
\pi_1^* \left( \Bfx[k\,;\mu] \right) \, $,  for all
%
% $ \, k= 1, 2, \dots, [f/2] \, $
%
 $ k $  and
 all
$ \, \mu \vdash \big(f\!-2k\big) \, $.

\saltino

  The next result gives us the required information about the splitting
of  $ \Sfx \, $.
%
% We prove it basing upon Brown's results, in a way
% independent of [HW1--2], and valid for  $ \, x = 0 \, $  too.
%

\salto

\proclaim{Proposition 3.2}  The algebra  $ \Sfx $  splits as  $ \;\; \Sfx
= \, \bigoplus\limits_{k=0}^{[f/2]} \bigoplus\limits_{\mu \vdash (f-2k)}
\Sfx[k\,;\mu] \; $,  \;\; where every  $ \, \Sfx[k\,;\mu] $  is a non-zero
simple algebra, unless  $ f $  is even and  $ \, (x,k) = \big(0,f/2\big)
\, $:  in this case,  $ \; {\Cal S}_f^{(0)}\big[f/2\,;(0)\big] \, \equiv
\, {\Cal S}_f^{(0)} \big[f/2\big] \, = \, 0 \; $.
\endproclaim

\demo{Proof}  First suppose  $ \, \Sfx[k\,;\mu] \neq 0 \, $.  As
$ \Bfx[k\,;\mu] $  is a generalized matrix algebra (Theorem 2.5), the
same is true for  $ \Sfx[k\,;\mu] $  too; but  $ \Sfx[k\,;\mu] $  is
semisimple, by construction, hence   --- \S 2.4(1) ---   it must be
simple.  Second, from the construction in \S 3.1 we get also
  $$  \displaylines{
   {} \indent   {\textstyle \bigoplus_{i \in J_k}} S_i  \; = \;
\Sfx[k]  \; = \;  \pi_1^* \left( \Bfx[k] \right)  \; = \;  \pi_1^*
\left( {\textstyle \bigoplus_{\mu \vdash (f-2k)}} \, \Bfx[k\,;\mu]
\right)  \; =   \hfill {}  \cr
   {} \hfill   = \;  {\textstyle \sum_{\mu \vdash (f-2k)}} \, \pi_1^*
\! \left( \Bfx[k\,;\mu] \right)  \; = \;  {\textstyle \sum_{\mu \vdash
(f-2k)}} \, \Sfx[k\,;\mu] \quad .   \indent {}  \cr }  $$
The summands in left hand sides are two-sided simple ideals, and the same is
true on right hand side: but the sum on the left is direct, and this easily
implies that each  $ \Sfx[k\,;\mu] $  is one of the  $ S_i $  and vice versa,
so that  $ \, \Sfx[k] = \bigoplus_{\mu \vdash (f-2k)} \Sfx[k\,;\mu] \; $.
Finally, since  $ \, \Sfx = \bigoplus_{k=0}^{[f/2]} \Sfx[k] \; $,  \, we
conclude that the splitting in the claim does hold.
                                                            \par
   Now we show that  $ \, \Sfx[k\,;\mu] \neq 0 \; $  for all  $ k $  and
$ \mu $  when  $ \, (x,k) \neq \big( 0, f/2 \big) \, $.  By definitions,
\S 2.4, Theorem 2.5 and Proposition 2.7, we have that
  $$  \Sfx[k\,;\mu]  = \big\{0\big\}  \,\iff\,
\Rad\, \big(\Bfx[k\,;\mu]\big) = \Bfx[k\,;\mu]  \,\iff\,
\text{\it rk}\Big(\Phi\big(\Bfx[k\,;\mu]\big)\Big) = 0  $$
and the last condition on the right clearly holds if and only if the
matrix  $ \, \Phi\big(\Bfx[k\,;\mu]\big) \, $  of all structure constants
of  $ \, \Bfx[k\,;\mu] \, $  is zero.  But this occurs exactly if and only
if  $ \, x = 0 \, $  and  $ \, k = f \big/ 2 \, $,  for even  $ f \, $.
In all other cases one has  $ \, \Sfx[k\,;\mu] \not= 0 \, $,  as
claimed.   $ \square $
\enddemo

\saltino

\proclaim{Corollary 3.3}  If  $ \, f \in \N \, $  is even, then
$ \; \text{\it Rad} \left( \Bfze \big(f/2\big) \right) = \, \Bfze
\big(f/2\big) = \, \Bfze \big[f/2\big] \; $.
\endproclaim

\saltino

  Indeed, the above also follows easily when remarking that
$ \, \Bfze\big[f/2\big] = \Bfze\big(f/2\big) \, $,  \, and
$ \Bfze\big(f/2\big) $  is just the  $ \Bbbk $--vector  space
$ \, \Bbbk^{{(f-1)!!}^2} \, $  endowed with the trivial
multiplication.

\salto

   {\bf 3.4 Semisimplicity of  $ \Bfx $.} \,  A general criterion for the
semisimplicity of  $ \Bfx $  is given in [Ru], [RS].  For the cases we shall
deal with, it reads as follows: if  $ \, x = n \in \N_+ \, $,  \, then
 \vskip3pt
   \centerline{\it  $ \Bfor $  and  $ \Bfsp $  are semisimple  $ \iff n \geq f-1
\; $  and  $ \, \text{\sl Char}\,(\Bbbk) = 0 \, $  or  $ \, \text{\sl Char}\,
(\Bbbk) > f \, $}
                                             \par
   \centerline{\it  $ \Bfze $  is semisimple  $ \iff f \in \{1,3,5\} \; $  and
$ \, \text{\sl Char}\,(\Bbbk) = 0 \, $  or  $ \, \text{\sl Char}\,(\Bbbk)
> f \, $}

\Salto

\centerline{ \bf  \S 4  Brauer algebras in Invariant Theory }

\ssalto

   {\bf 4.1 The Fundamental Theorems of Invariant Theory.} \
%
% We begin with some (more or less well-known) facts of Classical
% Invariant Theory; the general source is [We], but we shall also
% mention more specific   --- and more recent ---   references.
%                                                 \par
%
   Let  $ \, f \in \N_+ \, $  and  $ \, n \in \N \, $.  Let  $ V $
be a  $ \Bbbk $--vector  space of dimension  $ n $,  endowed with a
non-degenerate symmetric bilinear form  $ \, (\,\ , \ ) \, $,  and
let  $ O(V) $  be the associated orthogonal group.  Also, let  $ W $
be a  $ \Bbbk $--vector  space of dimension  $ 2n $,  endowed with a
non-degenerate skew-symmetric bilinear form  $ \, \langle \,\ , \ \rangle
\, $,  and let  $ Sp(W) $  be the associated symplectic group.  There exist
canonical
 isomor-\break{}phisms
 $ \, V {\buildrel \cong \over \longrightarrow} V^*,
v \mapsto (v,\; \cdot \; ) \, $,  $ \, W {\buildrel \cong \over \longrightarrow} W^*, w \mapsto \langle w, \; \cdot \; \rangle \, $,  \, which yield also
isomorphisms
  $$  \matrix
   \Theta_V : V \otimes V {\buildrel \cong \over \lllongrightarrow} End(V)
&  \Theta_W : W \otimes W {\buildrel \cong \over \lllongrightarrow} End(W)
\\
   \, v_1 \otimes v_2 \mapsto \Theta_V \left( v_1 \otimes v_2 \right)
\Big( v \mapsto \big( v_1, v \big) v_2 \Big) {\ }  &   {\ } w_1 \otimes
w_2 \mapsto \Theta_W \left( w_1 \otimes w_2 \right)  \Big( w \mapsto
\big\langle w_1 , w  \big\rangle w_2 \Big) \,  \\
      \endmatrix  $$
%
% Then  $ \, V^{\otimes 2f} {\buildrel \cong \over \longrightarrow} {\big(
% V^{\otimes 2f} \big)}^* $,  $ \, V^{\otimes 2f} = V^{\otimes f} \otimes
% V^{\otimes f} {\buildrel \cong \over \longrightarrow} End \big( V^{\otimes f}
% \big) $,  and  $ \, {\big( V^{\otimes 2f} \big)}^* {\buildrel \cong \over
% \longrightarrow} {End} \big( V^{\otimes f} \big) $,  whence also  $ \Psi_V :
% \Big( {\big( V^{\otimes 2f} \big)}^* \Big)^{O(V)} {\buildrel \cong \over
% \longrightarrow} \Big( {End} \big( V^{\otimes f} \big) \Big)^{O(V)} \! =
% \endor $;  and similarly for  $ W $, in particular  $ \, \Psi_W : \Big( {\big(
% W^{\otimes 2f} \big)}^* \Big)^{Sp(W)} {\buildrel \cong \over \longrightarrow}
% \Big( {End} \big( W^{\otimes f} \big) \Big)^{Sp(V)} \! = \endsp \, $.
%                                                  \par
%
   \indent   In this setting, we define  $ \; \psi_V := \Theta_V^{\,-1}
\left( {id}_V \right) \in V \otimes V \; $  and  $ \; \psi_W :=
\Theta_W^{\,-1} \left( {id}_W \right) \in W \otimes W \, $.

\salto
 \vskip3pt

\proclaim{Definition 4.2}  Fix  $ f \in \N_+ \, $.  For each pair  $ \,
p, q \in \{1, 2, \dots, f \} \, $  with  $ \, p \neq q \, $  we define
                                                            \par
   (a) \, a linear  {\sl contraction}  operator  $ \; \Phi_{p,q} :
V^{\otimes (f+2)} \! \llongrightarrow V^{\otimes f} \; $
(for  $ \, p < q \, $,  say), given by
 $ \;\; \Phi_{p,q} \big( v_1 \otimes v_2 \otimes \cdots \otimes v_{f+2} \big)
\, = \, \big( v_p \, , \,v_q \big) \cdot v_1 \otimes \cdots \widehat{v_p}
\otimes \cdots \otimes \widehat{v_q} \otimes \cdots \otimes v_{f+2} \;\; $;
 \eject
   (b) \, a linear  {\sl insertion}  operator  $ \; \Psi_{p,q} :
V^{\otimes f} \! \llongrightarrow V^{\otimes (f+2)} \; $,  obtained
by inserting the element  $ \psi_V $  in the positions  $ p $,
$ q \, $;
                                                            \par
   (c) \, an operator  $ \; \tau_{p,q} : V^{\otimes f} \llongrightarrow
V^{\otimes f} \; $  defined by  $ \; \tau_{p,q} := \Psi_{p,q} \circ
\Phi_{p,q} \; \big( \in {End} \big( V^{\otimes f} \big) \, \big) \, $.
 \vskip3pt
   The same definitions with  $ \langle \ , \ \rangle $  instead of
$ (\ ,\ ) $  give operators  $ \; \Phi_{p,q} : W^{\otimes (f+2)}
\longrightarrow W^{\otimes f} $,  $ \, \Psi_{p,q} : W^{\otimes f}
\longrightarrow W^{\otimes (f+2)} $,  $ \, \tau_{p,q} : W^{\otimes f}
\longrightarrow W^{\otimes f} \, $  in the symplectic case.
\endproclaim

\salto

   In addition, recall that the symmetric group  $ S_f $  acts on
$ V^{\otimes f} $  or  $ W^{\otimes f} $  by
 \vskip-5pt
  $$  \sigma \, : \, u_1 \otimes u_2 \otimes \cdots \otimes u_f \, \mapsto
\, u_{\sigma^{-1}(1)} \otimes u_{\sigma^{-1}(2)} \otimes \cdots \otimes
u_{\sigma^{-1}(f)}  \qquad  \forall\;\;  \sigma \in S_f  $$

\vskip3pt

   These constructions are connected with Brauer algebras, as we now explain.  The connection comes from a classical result over  $ \Bbb{C} $,  which later has been generalized to in [DP].  The technical condition required there,  {\sl for a given, fixed  $ \, f \in \N_+ \, $},  is the following:
  $$  \text{\sl Every polynomial  $ \, p(x) \in \Bbbk[x] \, $  of degree
$ f $  which vanishes on  $ \Bbbk $  is identically 0 \, .}   \eqno (4.1)  $$

\vskip-3pt

   {\bf In this section, we assume that the field  $ \Bbbk $  and
$ \, f \in \N_+ $  satisfy condition (4.1).}

\vskip7pt

  For instance, if  $ \, \text{\sl Char}\,(\Bbbk) = 0 \, $  or  $ \,
\text{\sl Char}\,(\Bbbk) > f \, $,  \, then  $ \Bbbk $  does satisfy
(4.1).  Actually, thanks to [DP], we can even assume  $ \Bbbk $  to
be any unital commutative ring satisfying (4.1).

\salto

   {\bf 4.3  Brauer algebras versus centralizer algebras.} \  When the
parameter  $ x $  is an integer, the Brauer algebra  $ \Bfx $  is strictly
related with the invariant theory for the orthogonal or the symplectic
groups.  Indeed, it is a ``lifting'' of one of the centralizer algebras
$ \endor $  or  $ \endsp \, $,  in the sense of the following result:

\salto

\proclaim{Theorem 4.4  {\rm (cf.~[Br], [DP])}}  Let  $ \, n \in \N_+ \, $,  and
let  $ V $  and  $ W $  respectively be an  $ n $--dimen\-sional  orthogonal
vector space and a  $ 2n $--dimensional  symplectic  vector space over
$ \, \Bbbk \, $.  Then there exist well-defined  $ \, \Bbbk $--algebra
epimorphisms,
 which are isomorphisms iff  $ \, n \geq f \, $,
  $$  \matrix
   \pi_V : \; \Bfor \llongtwoheadrightarrow \endor  \; \quad {}  &
{} \quad \;  \pi_W : \; \Bfsp \llongtwoheadrightarrow \endsp  \\
   \d_\sigma \mapsto \sigma  \, ,  \;\, \h_{p,q} \mapsto
\tau_{p,q}  \quad {\ }  &  {\ } \quad  \d_\sigma \mapsto sgn(\sigma)
\, \sigma  \, ,  \;\, \h_{p,q} \mapsto -\tau_{p,q}  \\
      \endmatrix  $$
   \indent   Moreover,  $ \, {End}_{\Bfor} \! \left( V^{\otimes f} \right)
= \big\langle O(V) \big\rangle \; $  and  $ \; {End}_{\Bfsp} \! \left(
V^{\otimes f} \right) = \big\langle Sp\,(W) \big\rangle \, $,  \, where
$ \, \big\langle X \big\rangle \, $  denotes the subalgebra of  $ \,
{End}_{\,\Bbbk} \! \left( U^{\otimes f} \right) $   ---  for  $ \,
U \in \big\{ V, W \big\} $  ---   generated by  $ X \, $.
\endproclaim

\vskip7pt

  The previous theorem   --- which follows from the First Fundamental Theorem
of Invariant Theory for  $ O(V) $  and  $ Sp(V) $  ---   concerns either  {\sl
positive\/}  or  {\sl even negative\/}  values of  $ x \, $.  The case of
{\sl odd negative\/}  parameter can be reduced to the odd positive case:
see [Wz], Corollary 3.5.  Finally, we shall cope with the case  $ \, x
= 0 \, $  through a direct approach.
                                              \par
%
%   The previous theorem (which roughly speaking follows from the First Fundamental Theorem of Invariant Theory for  $ O(V) $  and  $ Sp(V) \, $)  concerns either  {\sl positive\/}  or  {\sl even negative\/}  values of  $ x \, $.  The case of  {\sl odd
% negative\/}  parameter can be somehow reduced to the odd positive case:
% see [Wz], Corollary 3.5.  Finally, we shall cope with the case  $ \, x
% = 0 \, $  through a direct approach.  Therefore, our strategy now is to
% ``capitalize'' upon Theorem 4.4 above.
%                                               \par
%
   In order to describe the kernels of  $ \pi_V $  and  $ \pi_W \, $,  we
introduce some new objects.

\salto

   {\bf 4.5  Diagrammatic minors and diagrammatic Pfaffians.} \  Let us consider
the polynomial rings (in the symmetric or antisymmetric variables  $ x_{ij} $)
 \vskip-17pt
  $$  A^O := \Bbbk{[x_{ij}]}_{i,j=1,i \neq j}^{2f} \Big/ \left( x_{ij} =
x_{ji}  \right) \; ,  \qquad   A^{Sp} := \Bbbk{[x_{ij}]}_{i,j=1,i \neq j}^{2f}
\Big/ \left( x_{ij} = -x_{ji} \right)  $$
 \vskip-5pt
\noindent
 For  $ \, X \! \in \! \{O,{Sp}\} $,  define  $ A^X_f $  (the space of
{\sl multilinear elements}  in $ A^X $)  to be the  $ \Bbbk $--span  of all
monomials (of degree  $ f \, $)  $ \, x_{i_1 j_1} x_{i_2 j_2} \cdots x_{i_f
j_f} \, $  such that  $ (i_1, j_1, i_2, j_2, \dots, i_f, j_f) $  is a
permutation of  $ \, \big\{ 1, 2, 3, 4, \dots, 2f \,\big\} \, $.  Clearly,
$ A^X_f $  has a natural structure of  $ S_{2f} $--module.
 \vskip4pt
   In general, given two  $ f $--tuples  $ \, \text{\rm \bf i} := (i_1, i_2,
\dots, i_f) $  and  $ \text{\rm \bf j} := (j_1, j_2, \dots, j_f) $  such that
$ \, \{i_1, \dots, i_f\} \cup \{j_1, \dots, j_f\} = \{1,2,\dots,2f-1,2f\}
\, $,  \, we write  $ \, x_{\text{\, {\bf i}, {\bf j}}} := x_{i_1 j_1} x_{i_2 j_2}
\cdots x_{i_f j_f} \; $.  Then we define  $ \Bbbk $--vector  space isomorphisms
 $ \, \Phi_V \colon \, A_f^O {\buildrel \cong \over \longrightarrow} \, \Bfor \, $,
\; via  $ \, x_{\text{\, {\bf i}, {\bf j}}} \mapsto \d_{\text{{\bf i}, {\bf j}}} \, $,
and
 $ \, \Phi_W \colon \, A_f^{Sp} {\buildrel \cong \over \longrightarrow} \, \Bfsp \, $,
\; via  $ \, x_{\text{\, {\bf i}, {\bf j}}} \mapsto \varepsilon(\d_{\text{{\bf i}, {\bf j}}}) \cdot \d_{\text{{\bf i}, {\bf j}}} \, $,  \;
%
%   $$  \matrix
%    {} \qquad  \Phi_V : A_f^O {\buildrel \cong \over \llongrightarrow}
% \Bfor  \qquad {}  &  {} \qquad  \Phi_W : A_f^{Sp} {\buildrel \cong \over
% \llongrightarrow} \Bfsp  \qquad {}  \\
%    x_{\text{\, {\bf i}, {\bf j}}} \mapsto \d_{\text{{\bf i}, {\bf j}}}  &
% x_{\text{\, {\bf i}, {\bf j}}} \mapsto \varepsilon(\d_{\text{{\bf i},
% {\bf j}}}) \cdot \d_{\text{{\bf i}, {\bf j}}}  \\
%       \endmatrix  $$
%
where  $ \varepsilon(\d_{\text{{\bf i},{\bf j}}}) $  is the ``sign'' of
$ \d_{\text{{\bf i},{\bf j}}} $  defined as in \S 1.4.  Using them, an
$ S_{2f} $--action  is defined on  $ \Bfor $,  resp.~$ \Bfsp $,  based
upon that on  $ A_f^O $,  resp.~$ A_f^{Sp} $,  letting  $ \, \sigma \in
S_{2f} \, $  act on  $ \Bfor $,  resp.~$ \Bfsp $,  as  $ \, \Phi_V \circ
\sigma \circ \Phi_V^{-1} \, $,  resp.~$ \, \Phi_W \circ \sigma \circ
\Phi_W^{-1} \, $ (this action is studied in depth in [Hu] and in [DH]).

\salto

\noindent   {\bf Definition 4.6}
                                        \par
  {\it (a) \, We call {\bf (diagrammatic) minor}  of order  $ r $
($ \in \N_+ $)  every element of  $ \, \Bfx $  which is the image
through  $ \Phi_V $  of an element of type
  $$  {\textstyle \sum_{\sigma \in S_r}} \, sgn(\sigma) \cdot x_{i_1 j_{\sigma(1)}}
x_{i_2  j_{\sigma(2)}} \cdots x_{i_r j_{\sigma(r)}} \cdot x_{i_{r+1}
j_{r+1}} x_{i_{r+2} j_{r+2}} \cdots x_{i_{f-1} j_{f-1}} x_{i_f j_f}
\hskip11pt   \eqno (4.2)  $$
with  $ \; \{i_1, \dots, i_r\} \cup \{j_1, \dots, j_r\}
\cup \{i_{r+1}, \dots, i_f\} \cup \{j_{r+1}, \dots, j_f\}
= \{1, 2, 3, \dots, 2f\} \; $.
                                        \par
   We denote by  $ \, \text{\sl Min}^{\,(x)}_{f;\,r} $  the set of all
(diagrammatic) minors  of order  $ r $  in  $ \, \Bfx \, $.
                                        \par
  (b) \, We call {\bf (diagrammatic) Pfaffian}  of order  $ 2r $  ($ \in 2
\N_+ $)  every element of  $ \Bfx $  which is the image through  $ \Phi_W $
of an element of type
  $$  {\textstyle \sum_{\Sb
              h_1 < k_1, h_2 < k_2, \dots  \\
              h_1 < h_2 < h_3 < \cdots
            \endSb}} \hskip1pt sgn \hskip2pt
    {\textstyle
  \Big(\hskip-4pt
        {{\hskip5pt 1 \hskip9pt 2 \hskip9pt \ldots \hskip8pt {2r-1} \hskip5pt 2r}
                     \atop
        {\hskip5pt h_1 \hskip5pt k_1 \hskip10pt \ldots \hskip9pt h_r \hskip9pt k_r}} \Big)} \cdot
%
%     \pmatrix
%       1 \!\!  &  \!\! 2 \! &  \! \ldots \!  &  \! 2r \! - \! 1 \!
% &  \! 2r  \\
%       h_1 \!\!  &  \!\! k_1 \!  &  \! \ldots \! &  \! h_r \!
% &  \! k_r
%     \endpmatrix \,
%
x_{h_1 k_1} x_{h_2 k_2} \cdots x_{h_r k_r} x_{i_{r+1}
j_{r+1}} \cdots x_{i_f j_f}   \eqno (4.3)  $$
with  $ \, \{h_1, \dots, h_r\} \cup \{k_1, \dots, k_r\}
\cup \{i_{r+1}, \dots, i_f\} \cup \{j_{r+1}, \dots, j_f\}
= \{1, 2, 3, \dots, 2f\} \, $.
                                        \par
   We denote by  $ \, \text{\sl Pf}^{\;(x)}_{f;\,r} $  the set of all
(diagrammatic) Pfaffians  of order  $ \, 2 \, r $  in  $ \, \Bfx \, $.
                                        \par
  (c) \, If  $ X $  is any given (diagrammatic) minor or Pfaffian,
we call  {\sl fixed edge}  of  $ X $  any edge which occurs the same in
all diagrams occurring in the expansion of  $ X \, $.  We call  {\sl fixed
vertex}  of  $ X $  any vertex (in  $ \V_f $)  belonging to a fixed edge
of  $ X \, $.  We call  {\sl fixed part}  of  $ X $  the datum of all fixed
edges and all fixed vertices of  $ X \, $.
                                        \par
  (d) \, If  $ X $  is any given (diagrammatic) minor or Pfaffian, we call
{\sl moving vertex}  of  $ X $  any vertex (in  $ \V_f $)  which is not
fixed in  $ X $.  We call  {\sl moving part}  of $ X $  the datum of all
vertices which are not fixed in  $ X $.}

\salto

   {\bf Remarks 4.7.} \; {\it (a)} \, By definitions and Proposition
4.6, any diagrammatic minor of order  $ r $  is an  {\sl alternating sum
of  $ f $--diagrams\/}:  to be precise, it is an  $ S_r $--antisymmetric
sum of  $ f $--diagrams.  On the other hand, due to the sign entering in the
definition of  $ \alpha_W \, $,  all diagrams occurring in the expansion
of a diagrammatic Pfaffian appear there with the same sign.  Thus,
each diagrammatic Pfaffian is (up to sign) just a simple
{\sl sum of  $ f $--diagrams}.
                                                   \par
   {\it (b)} \, Let  $ \delta_r $  be a minor of order  $ r \, $.  Its
moving vertices may be partitioned into two sets  $ I $,  $ J $  (each
of  $ r $  elements) so that, looking at all the diagrams occurring in
the expansion of  $ \delta_r \, $,  no vertex in one of these sets is
ever joined to a vertex in the same set, but it is joined to each
vertex in the other set.  These  $ I $  and  $ J $  correspond, via
$ \Phi_V $,  to the set of rows and the set of columns (or vice versa) in
the matrix  $ {\big( x_{ij} \big)}_{i,j=1}^{2f} $  on which the minor
corresponding to  $ \delta_r $  is computed.  So in the sequel
expressions like  ``$ v $  is a row vertex and  $ w $  is a column
vertex'' will mean that  $ v $  and  $ w $  are moving
vertices which belong one to  $ I $  and the other to  $ J \, $.  Similarly,
by  ``$ v $  and  $ w $  are both row vertices'' or ``column vertices''
well mean that they are moving vertices which both belong to  $ I $  or
both to  $ J $.  Indeed, a minor  $ \delta_r $  is determined,
up to sign, by:  {\it (i)\/}  assigning its fixed part;  {\it (ii)\/}
assigning the sets  $ I $  and  $ J \, $,  each endowed with a labelling
of its vertices by  $ \{1, 2, \dots, r\} $;   {\it (iii)\/}  joining every
vertex in one set   --- say  $ I $ ---   to a vertex in the other set   ---
say  $ J $ ---   according to a permutation  $ \, \sigma \in S_r \, $,  so
to get an  $ f $--diagram  $ \d(\sigma) $;  {\it (iv)\/}  adding up the
diagrams  $ \d(\sigma) $  with coefficient  $ sgn(\sigma) $,  for all
$ \sigma \in S_r \, $:  this eventually gives  $ \pm \d_r $  (the sign
depends on the labellings).
                                                   \par
   {\it (c)} \, The step  {\it (iii)\/}  above may be better understood
as follows.  First, join every vertex in  $ I $  with the vertex in  $ J $
labelled with the same number: this gives the diagram  $ \d(id) $  which,
outside the fixed part, is given by the  $ r $  edges  $ \{i_1,j_1\},
\dots, \{i_r,j_r\} $  (where  $ \, \{i_1, \dots, i_r\} = I \, $,  $ \,
\{j_1, \dots, j_r\} = J \, $).  Second, let  $ S_r $  act on  $ J $,  and
let  $ \d[\sigma] $  be the diagram which is equal to  $ \d(id) $  in the
fixed part and outside it is given by the  $ r $  edges  $ \big\{i_1,
\sigma(j_1) \big\}, \dots, \big\{ i_r, \sigma(j_r) \big\} $:  then  $ \,
\d[\sigma] = \d(\sigma) \, $.  Then we can also write  $ \delta_r $  as
an  $ S_r $--antisymmetric  sum
  $$  \delta_r \; = \; {\textstyle \sum_{\sigma \in S_r}} sgn(\sigma)
\, \d(\sigma) \; = \; {\textstyle \sum_{\sigma \in S_r}} sgn(\sigma)
\, \d[\sigma] \; = \; {\textstyle \sum_{\sigma \in S_r}} sgn(\sigma)
\, \sigma.{\hskip0,5pt}\d[id]  $$
\indent
   {\it (d)} \, The counterpart for Pfaffians of  {\it (b)}  and  {\it (c)}
above is that every Pfaffian of order  $ 2r $  is the sum of all diagrams
obtained by assigning the fixed part and joining the  $ 2r $  vertices in
the moving part with  $ r $  edges in all possible ways.
                                        \par
  {\it (e)} \, Examples of diagrammatic minors or Pfaffians can be found
in [Ga], Example 3.6.

\salto

   The importance of diagrammatic minors and Pfaffians lies in the following:
%
% reformulation of the Second Fundamental Theorem of Invariant Theory
%

\ssalto

\proclaim{Theorem 4.8}  {\sl ([Ga], \S 3)}  Let  $ \, n \in \N_+ \, $  and
$ \, f \in \N_+ \, $  with (4.1) satisfied by the field  $ \, \Bbbk \, $.
                                               \par
   (a) \, The kernel of  $ \, \pi_V : \Bfor \! \llongtwoheadrightarrow
\endor \, $  is the  $ \, \Bbbk $--span  of the set of all diagrammatic
minors in  $ \Bfor $  of order  $ \, n+1 \, $.  In particular, it is an
$ S_{2f} $--submodule.
                                               \par
   (b) \, The kernel of  $ \, \pi_W \colon \Bfsp \! \llongtwoheadrightarrow
\endsp \, $  is the  $ \, \Bbbk $--span  of the set of all diagrammatic
Pfaffians in  $ \Bfsp $  of order  $ \, 2(n+1) \, $.  In particular, it
is an  $ S_{2f} $--submodule.
\endproclaim

\saltino

   {\bf 4.9 Comparison with others' work.} The problem of describing  $ \Ker\,
(\pi_W\!) $  is solved also by Hu in [Hu].  He studies the action of 
$ S_{2f} $  onto  $ D_f \, $,  and the  $ S_{2f} $--module  structure
of the  $ \Z $--algebra  $ \Bfx(\Z) \, $,  \, the  $ \Z $--span
of  $ D_f \, $.  As main results, he finds a new basis, and Specht filtrations, for  $ \Bfx(\Z) \, $,  and a characteristic free description of  $ \Ker\,(\pi_W) $,  proving that it is an  $ S_{2f} $--submodule  of  $ \Bfx(\Z) \, $.  Some of his results can be compared to ours: for instance, Theorem 3.4 in [Hu], describing 
$ \Ker\,(\pi_W) $,  coincides with our Theorem 4.8{\it (b)}.
                                               \par
   Another interesting point is Lemma 3.3 in [Hu], which proves that certain
sums of diagrams do belong to  $ \Ker\,(\pi_W) \, $.  Well, definitions imply
that any such sum is simply (a special type of) a  {\sl Pfaffian}   --- of
order 2(a+b) ---   in the sense of our  Definition 4.6{\it (b)\/}  (see
Remarks 4.7 too).  Therefore, Hu's lemma is just a special case of our
Theorem 4.8{\it (b)}.
                                               \par
   Similarly, an analogous solution for  $ \Ker\,(\pi_V) $  is (just very recently) provided in [DH].   
%  
% \Salto   
%  
 \eject   

\centerline{ \bf  \S 5  Within the radical of  $ \, \Bfx $ }

\ssalto

   {\bf 5.1 From Invariant Theory to  $ \Rad \left( \Bfx \right) \, $.}
\  In the present work, relying on the results of Invariant Theory in
\S\S 3--4, we locate a large family of elements in  $ \Rad \left( \Bfx
\right) $,  namely diagrammatic minors or Pfaffians,  when  $ x $  is
an integer which is not odd negative.

\vskip7pt

   Until \S 5.15,  {\bf we assume now  $ \text{\sl Char}\,(\Bbbk) = 0
\, $}.  {\sl In particular, (4.1) holds for any  $ f \in \N_+ \, $}.

\vskip7pt

   Being in characteristic zero, the orthogonal groups are linearly reductive.
Hence, by general theory,  $ \endor $  is semisimple; so the epimorphism
$ \, \pi_V \, \colon \Bfor \!\! \relbar\joinrel\twoheadrightarrow
\! \endor $  factors through  $ \, \pi_1 \, \colon \, \Bfor \!\!
\relbar\joinrel\relbar\joinrel\twoheadrightarrow \Sfor \, $.  In other
words,  $ \pi_V $  is the composition of maps  $ \, \pi_V = \pi_2 \circ
\pi_1 \, \colon \, \Bfor \!\! \relbar\joinrel\twoheadrightarrow \Sfor \!\!
\relbar\joinrel\twoheadrightarrow \endor \, $  \, where  $ \pi_2 $  is the
map given by the universality of the semisimple quotient.  It follows that
$ \, \Rad \left( \Bfor \right) = \Ker\, \big( \pi_1 \big) \subseteq \Ker\,
\big( \pi_V \big) \, $.  By  Theorem 4.8{\it (a)},  the latter space is the
$ \Bbbk $--span  of all (diagrammatic) minors of order  $ n+1 $  in  $ \,
\Bfor $.  So we shall look here for elements of  $ \, \Rad \left(\! \Bfor
\right) $,  in particular we shall determine (in Theorem 5.3 and Theorem 5.5)
exactly which of those minors do belong to  $ \Rad \left(\! \Bfor \right) $.
                                                \par
   The same arguments give also  $ \, \pi_W = \pi_2 \circ \pi_1 \, \colon
\, \Bfsp \!\! \relbar\joinrel\relbar\joinrel\twoheadrightarrow \Sfsp \!\!
\relbar\joinrel\relbar\joinrel\twoheadrightarrow \endsp \, $,  \, so  $ \,
\Rad \left( \Bfsp \right) = \Ker\, \big( \pi_1 \big) \subseteq \Ker\,
\big( \pi_W \big) \, $.  The latter space is the  $ \Bbbk $--span  of all
(diagrammatic) Pfaffians of order  $ 2(n+1) $  in  $ \, \Bfsp \, $,  by
Theorem  4.8{\it (b)}.  Hence we shall look here for elements of  $ \,
\Rad \left( \Bfsp \right) \, $,  \, and we shall determine exactly   --- in
Theorem 5.3 and Theorem 5.5 again   ---   which of those Pfaffians
actually do belong to  $ \, \Rad \left( \Bfsp \right) \, $.

\salto
 \vskip3pt

\proclaim{Proposition 5.2}  Let  $ \, n \in \N_+ \, $,
$ \, f \in \N_+ \, $.  Then
  $$  \leqalignno{
   \Ker\,(\pi_V) \; {\textstyle \bigcap} \; \Big( {\textstyle
\bigoplus_{h=0}^{[f/2]}} \, {\textstyle \bigoplus_{\Sb
         \mu \vdash (f-2h)  \\
         \mu_1^t + \mu_2^t \leq n  \\
           \endSb}} \Bfor \langle h \, ; \mu \rangle \Big)
\;  &  \subseteq \; \Rad\, \Big( \Bfor \Big)  &  \quad (a)  \cr
   \Bfze \big( \big[ (f+1)/2 \big] \big) \; \subseteq \;
\Rad\, \Big( \Bfze \Big)  &   &  \quad (b)  \cr
   \Ker\,(\pi_W) \; {\textstyle \bigcap} \; \Big( {\textstyle
\bigoplus_{h=0}^{[f/2]}} \, {\textstyle \bigoplus_{\Sb
         \mu \vdash (f-2h)  \\
         \mu_1^t \leq n  \\
           \endSb}} \Bfsp \langle h \, ; \mu \rangle \Big)
\;  &  \subseteq \; \Rad\, \Big( \Bfsp \Big)  &  \quad (c)  \cr }  $$
\endproclaim

\demo{Proof} {\it (a)} \,  By [Wz], \S 3, we know that  $ \endor $  splits
into a direct sum of pairwise non-isomorphic simple algebras as  $ \, \endor
= \bigoplus_{h=0}^{[f/2]} \bigoplus_{\Sb
             \mu \vdash (f-2h)  \\
             \mu_1^t + \mu_2^t \leq n  \\
           \endSb} A[h\,;\mu] \, $,  \, where  $ \mu^t $  is the dual partition
to  $ \mu $,  as in \S 2.2.  Furthermore, the analysis in \S\S 3--4 shows that
$ \; A[h\,;\mu] \, = \, \pi_2 \left( \Sfor[h\,;\mu] \right) \, = \, \pi_2 \left( \pi_1^* \left( \Bfor[h\,;\mu] \right) \right) \; $.
%
%   $$  A[h\,;\mu] \, = \, \pi_2 \left( \Sfor[h\,;\mu] \right)
% \, = \, \pi_1^* \left( \Bfor[h\,;\mu] \right) \;\; .  $$
%
Now, the map
  $$  \pi_2 \;\; \colon \; {\textstyle \bigoplus\limits_{h=0}^{[f/2]}}
{\textstyle \bigoplus\limits_{\mu \vdash (f-2h)}} \! \Sfor[h\,;\mu] \;
= \; \Sfor \relbar\joinrel\relbar\joinrel\twoheadrightarrow \;
\endor \; = \; {\textstyle \bigoplus\limits_{h=0}^{[f/2]}}
{\textstyle \bigoplus_{\Sb
             \mu \vdash (f-2h)  \\
             \mu_1^t + \mu_2^t \leq n  \\
           \endSb}}  \, A[h\,;\mu]  $$
 \eject
\noindent
 (cf.~Proposition 3.2) must have kernel  $ \; \Ker\,(\pi_2)
\, = \, \bigoplus\limits_{h=0}^{[f/2]} \bigoplus_{\Sb
         \mu \vdash (f-2h)  \\
         \mu_1^t + \mu_2^t > n  \\
           \endSb}  \! \Sfor[h\,;\mu] \; $,  \; and it must map
$ \bigoplus\limits_{h=0}^{[f/2]} \! \bigoplus_{\Sb
         \!\! \mu \vdash (f-2h)  \\
         \!\! \mu_1^t + \mu_2^t \leq n  \\
           \endSb}  \!\! \Sfor[h\,;\mu] \, $  isomorphically onto
$ \bigoplus\limits_{h=0}^{[f/2]} \! \bigoplus_{\Sb
         \!\! \mu \vdash (f-2h)  \\
         \!\! \mu_1^t + \mu_2^t \leq n  \\
           \endSb}  \!\! A[h\,;\mu] = \endor $.
                                                \par
   Let now consider an element  $ \; y \in \Ker\,(\pi_V) \, \bigcap \,
\Bigg( \bigoplus_{h=0}^{[f/2]} \bigoplus_{\Sb
         \mu \vdash (f-2h)  \\
         \mu_1^t + \mu_2^t \leq n  \\
           \endSb} \, \Bfor \langle h \, ; \mu
\rangle \Bigg) \; $.  Then  $ \pi_1(y) $  belongs to
$ \, \bigoplus_{h=0}^{[f/2]}
\, \bigoplus_{\Sb
         \mu \vdash (f-2h)  \\
         \mu_1^t + \mu_2^t \leq n  \\
           \endSb} \Sfor[h\,;\mu] \; $.
But on the latter space  $ \pi_2 $  acts injectively, hence  $ \;
\pi_2\big(\pi_1(y)\big) = \pi_V(y) = 0 \; $  implies  $ \; \pi_1(y)
= 0 \; $,  \; so  $ \, y \in \Ker\, \big( \pi_1 \big) \equiv
\Rad \left( \Bfor \right) \, $,  q.e.d.

\vskip7pt

   {\it (b)} \, If  $ f $  is odd the claim is empty, and there is nothing
to prove.  If  $ f $  is even, then  $ \, \big[ (f+1) / 2 \big] = f/2 \, $,
and the claim follows from Corollary 3.3.  Indeed, the latter gives  $ \,
\Bfze\big[f/2\big] \equiv \Bfze\big(f/2\big) = \Rad\, \Big( \Bfze\big(
f/2 \big) \Big) \, \subseteq \Rad \left( \Bfze \right) \, $,  with the
last inclusion following by easy arguments of Artinian algebras (as in
[HW1], \S 4.B).  Otherwise, one can proceed as follows.  Definitions imply
$ \, \Bfze \big( f/2 \big) \subseteq \text{\it Ann}\, \big( H^\mu_{f,h}
\big) \, $  for all  $ h $,  $ \mu \vdash (f-2h) \, $;  therefore  $ \Bfze
\big( f/2 \big) $  kills also all simple  $ \Bfze $--modules,  for they
are quotients of the  $ H^\mu_{f,h} $'s.  But the radical of a finite
dimensional  $ \Bbbk $--algebra  $ A $  is characterized (or defined)
by  $ \; \Rad\,(A) = \bigcap_{M \in Spec(A)} Ann(M) \; $,  \; where
$ Spec(A) $  is the set of finite dimensional simple  $ A $--modules.
Thus we conclude that  $ \, \Bfze\big(f/2\big) \subseteq \Rad \left(
\Bfze \right) \, $,  whence the claim follows again.

\vskip7pt

   {\it (c)} \, By [Wz], \S 3, there is also a splitting  $ \; \endsp
= \, \bigoplus_{h=0}^{[f/2]} \bigoplus_{\Sb
                           \mu \vdash (f-2h)  \\
                           \mu_1^t \leq n  \\
  \endSb} A[h\,;\mu] \; $,  \; and the analysis in \S\S 3--4 ensures that
$ \; A[h\,;\mu] \, = \, \pi_2 \left( \Sfsp[h\,;\mu] \right) \, = \,
\pi_1^* \left( \Bfsp[h\,;\mu] \right) \; $.
%
% %
%  \vskip-7pt
% %
%   $$  A[h\,;\mu] \, = \, \pi_2 \left( \Sfsp[h\,;\mu] \right) \, = \,
% \pi_1^* \left( \Bfsp[h\,;\mu] \right) \; .  $$
% %
%
Like before, using the splitting of  $ \Cal{S}_f^{(x)} $  in Proposition 3.2,
we see that the map
 \vskip-13pt
  $$  \pi_2 \;\; \colon \; {\textstyle \bigoplus\limits_{h=0}^{[f/2]}}
{\textstyle \bigoplus\limits_{\mu \vdash (f-2h)}} \! \Sfsp[h\,;\mu]
\, = \,  \Sfsp \! \relbar\joinrel\relbar\joinrel\twoheadrightarrow
\, \endsp \, = \, {\textstyle \bigoplus\limits_{h=0}^{[f/2]}} \,
{\textstyle \bigoplus_{\Sb
             \! \mu \vdash (f-2h)  \\
             \! \mu_1^t \leq n  \\
           \endSb}}  A[h\,;\mu]  $$
 \vskip-7pt
\noindent
 must have kernel  $ \; \Ker\,(\pi_2) \, = \,
\bigoplus\limits_{h=0}^{[f/2]} \bigoplus_{\Sb
         \mu \vdash (f-2h)  \\
         \mu_1^t > n  \\
           \endSb}  \! \Sfsp[h\,;\mu] \; $;  \; in addition, it must map
$ \; \bigoplus\limits_{h=0}^{[f/2]} \! \bigoplus_{\!\Sb
         \mu \vdash (f-2k)  \\
         \mu_1^t \leq n  \\
           \endSb}  \! \Sfsp[h\,;\mu] \; $  isomorphically onto
$ \; \bigoplus\limits_{h=0}^{[f/2]} \! \bigoplus_{\!\Sb
         \mu \vdash (f-2h)  \\
         \mu_1^t \leq n  \\
           \endSb}  \! A[h\,;\mu] = \endor \, $.
                                                \par
   Let now  $ \; \eta \in \Ker\,(\pi_W) \, \bigcap \,
\Big( \bigoplus_{h=0}^{[f/2]} \, \bigoplus_{\Sb
         \mu \vdash (f-2h)  \\
         \mu_1^t \leq n  \\
           \endSb} \Bfsp \langle h \, ; \mu \rangle \Big) \; $.  Then
$ \pi_1(\eta) $  belongs to  $ \, \bigoplus\limits_{h=0}^{[f/2]} \,
\bigoplus_{\Sb
         \mu \vdash (f-2h)  \\
         \mu_1 \leq n  \\
           \endSb} \Sfsp[h\,;\mu] \; $.
As  $ \pi_2 $  acts injectively on the latter space, we get that  $ \;
\pi_2 \big( \pi_1 ( \varpi_{n+1} ) \big) = \pi_W ( \varpi_{n+1} ) = 0
\; $  yields  $ \; \pi_1(\eta) = 0 \; $,  \; so  $ \, \eta \in \Ker\,
\big( \pi_1 \big) \equiv \Rad \! \left( \Bfsp \right) \, $.
\hskip5pt  $ \square $
\enddemo

\saltino

   The previous statement has the following direct consequence:

\ssalto

\proclaim{Theorem 5.3}  Let  $ \, n \in \N_+ \, $,  $ \, f \in \N_+ \, $,
$ \, k := \left[ {{f-n+1} \over {2}} \right] \, $.  Then (notation of
Definition 4.6)
                                                 \par
   (a) \, every minor of order  $ (n+1) $  in  $ \Bfor (k) $
belongs to  $ \Rad \left( \Bfor \right) \, $,  \, hence
\vskip-3pt
  $$  \text{$ \Bbbk $--{\sl span}}\Big( \text{\sl Min}^{\,(n)}_{f;\,n+1}
\,{\textstyle \bigcap}\, \Bfor(k) \Big) \, \subseteq \, \Rad \left(
\Bfor \right)  \quad \text{;}  $$
\vskip-3pt
   \indent   (b) \, every  $ f $--diagram  in  $ \Bfze \! \left( \big[
(f+1)/2 \big] \right) \, $  belongs to  $ \Rad \left( \Bfze \right) \, $,
\, hence
\vskip-3pt
  $$  \Bfze \! \left( \big[ (f+1)/2 \big] \right) \,
\subseteq \, \Rad \left( \Bfze \right)  \quad \text{;}  $$
\vskip-3pt
   \indent   (c) \, every Pfaffian of order  $ 2(n+1) $  in  $ \Bfsp (k) $
belongs to  $ \Rad \left( \Bfsp \right) \, $,  \, hence
\vskip-3pt
  $$  \text{$ \Bbbk $--{\sl span}}\Big( \text{\sl Pf}^{\;(-2n)}_{f;\,n+1}
\,{\textstyle \bigcap}\, \Bfsp(k) \Big) \, \subseteq \, \Rad \left(
\Bfsp \right)  \quad \text{.}  $$
\vskip-3pt
\endproclaim

\demo{Proof} {\it (a)} \,  Let  $ \delta_{n+1} $  be a minor of order
$ (n+1) $  in  $ \Bfor(k) \, $,  with  $ \, k := \left[ {{f-n+1} \over {2}}
\right] \, $.  This means that all the  $ f $--diagrams occurring (with
non-zero coefficient) in the expansion of  $ \delta_{n+1} $  have at least
$ k $  arcs, so they have at most  $ f-2k $  vertical edges, with  $ \,
f-2k < n+1 \, $.  Therefore  $ \; \delta_{n+1} \in \bigoplus_{h=0}^{[f/2]}
\, \bigoplus_{\Sb
         \mu \vdash (f-2h)  \\
         \mu_1^t + \mu_2^t \leq n  \\
           \endSb} \Bfor \langle h \, ; \mu \rangle \; $,  \; and in
addition  $ \, \delta_{n+1} \in \Ker\,(\pi_V) \, $,  by  Theorem
4.8{\it (a)}.  Then  Proposition 5.2{\it (a)\/}  gives  $ \, \delta_{n+1}
\in \Rad \left( \Bfor \right) \, $,  \, as claimed.
                                                 \par
   {\it (b)} \,  This is obvious from  Proposition 5.2{\it (b)}.
                                                 \par
   {\it (c)} \,  Let  $ \varpi_{n+1} $  be a Pfaffian of order
$ 2(n+1) $  in  $ \Bfsp(k) \, $,  with  $ \, k := \left[ {{f-n+1} \over {2}}
\right] \, $.  Then, as in  {\it (a)},  all the  $ f $--diagrams occurring
(with non-zero coefficient) in the expansion of  $ \varpi_{n+1} $  have at
most  $ f-2k $  vertical edges, with  $ \, f-2k < n+1 \, $.  It follows
that  $ \; \varpi_{n+1} \in \bigoplus_{h=0}^{[f/2]}
\, \bigoplus_{\Sb
         \mu \vdash (f-2h)  \\
         \mu_1^t \leq n  \\
           \endSb} \Bfsp \langle h \, ; \mu \rangle \; $,  \; and moreover
$ \, \varpi_{n+1} \in \Ker\,(\pi_W) \, $,  by  Theorem 4.8{\it (b)}.
Thus  Proposition 5.2{\it (c)\/}  gives  $ \, \varpi_{n+1} \in \Rad
\left( \Bfsp \right) \, $,  \, as expected.   $ \square $
\enddemo

\salto

   Note that in the statement above, part  {\it (b)\/} is an improvement of
Corollary 3.3.  Moreover, it can be formulated as in  {\it (a)\/}  or in
{\it (c)},  with  $ \, n = 0 \, $;  \, indeed, a diagram is just a minor
of order  $ 1 $,  or a Pfaffian of order  $ 2 $,  and vice versa.
 \vskip3pt
   For the next step, we need a technical (combinatorial) result about
minors and Pfaffians:

\salto

\proclaim{Lemma 5.4}  Let  $ \, n \in \N_+ \, $,  $ \, f \in \N_+ \, $,  and
let  $ \Bfor $,  $ \Bfsp $  be defined over  {\sl any}  ring  $ \Bbbk \, $.
                                                        \par
   (a) \, Let  $ \,\d $  be an  $ f $--diagram,  and  $ \, \delta_{n+1} $
a minor of order  $ n+1 $  in  $ \Bfor \, $.  If  $ \, \d $  has an arc \
$ \; r^- \, \epsfbox{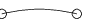} \, s^- \; $,  \  resp.~\  $ \; r^+
\, \epsfbox{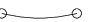} \, s^+ \; $,  \  and  $ r^+ $  and  $ s^+ $,
resp.~$ r^- $  and  $ s^- $,  are moving vertices in  $ \delta_{n+1} \, $,
then  $ \; \d \cdot \delta_{n+1} = 0 \; $,  \, resp.~$ \; \delta_{n+1} \cdot
\d = 0 \; $.  Otherwise,  $ \, \d \cdot \delta_{n+1} \, $,  resp.~$ \,
\delta_{n+1} \cdot \d \, $,  \, is a power of  $ n $  times a minor of
order  $ n+1 \, $.
                                                        \par
   Similarly, if  $ \, j \in J_{f,k} \, $  is an  $ (f,k) $--junction
(for some  $ k $)  having an arc \  $ \; r \, \epsfbox{r--s_bardown.eps}
\, s \; $,  \ and  $ r^- $  and  $ s^- $  are moving vertices in
$ \delta_{n+1} \, $,  then  $ \; \delta_{n+1}.j = 0 \; $  in
$ H^\mu_{f,k} $  for all  $ \, \mu \vdash (f-2k) \, $.
                                                        \par
   (b) \, Let  $ \,\d $  be an  $ f $--diagram,  and  $ \, \varpi_{n+1} $
a Pfaffian of order  $ 2(n+1) $  in  $ \Bfsp \, $.  If  $ \, \d $  has an
arc \  $ \; r^- \, \epsfbox{r--s_bardown.eps} \, s^- \; $,  \ resp.~\  $ \;
r^+ \, \epsfbox{r--s_barup.eps} \, s^+ \; $,  \ and  $ r^+ $  and  $ s^+ $,
resp.~$ r^- $  and  $ s^- $,  are moving vertices in  $ \varpi_{n+1} \, $,
then  $ \; \d \cdot \varpi_{n+1} = 0 \; $,  \, resp.~$ \; \varpi_{n+1} \cdot
\d = 0 \; $.  Otherwise,  $ \, \d \cdot \varpi_{n+1} \, $,  \, resp.~$ \,
\varpi_{n+1} \cdot \d \, $,  \, is a power of  $ (-2n) $  times a Pfaffian
of order  $ 2(n+1) \, $.
                                                        \par
   Similarly, if  $ \, j \in J_{f,k} \, $  is an  $ (f,k) $--junction  (for
some  $ k $)  having an arc \  $ \; r \, \epsfbox{r--s_bardown.eps} \, s \; $,
\  and  $ r^- $  and  $ s^- $  are moving vertices in  $ \varpi_{n+1} \, $,
then  $ \; \varpi_{n+1}.j = 0 \; $  in  $ H^\mu_{f,k} $  for all  $ \, \mu
\vdash (f-2k) \, $.
\endproclaim

\demo{Proof} Lemma 3.9 in [Ga] proves almost all of the present claim.  What
is missing is only the parts which start with ``Otherwise''.  Now, checking
also these facts is immediate from definitions.  Here we just mention
explicitly that a coefficient  $ n \, $,  or  $ (-2n) \, $,  will pop up
--- so its exponent will increase to give a power of  $ n \, $,  or  $ (-2n) $
---   whenever an arc in  $ \d $  matches a  {\sl fixed\/}  arc in the minor
$ \delta_{n+1} $   --- for  {\it (a)}  ---   or in the Pfaffian
$ \varpi_{n+1} $   --- for  {\it (b)}\,.   $ \square $
\enddemo

\salto

   Theorem 5.3 gives a sufficient condition for a minor, resp.~a Pfaffian
(of the proper order) to belong to  $ \Rad\left(\Bfor\right) \, $,  \,
resp.~to  $ \Rad\left(\Bfsp\right) \, $.  The next result claims that
this is necessary too.  Note that again claim  {\it (b)\/}  may be
enclosed in  {\it (a)\/}  or in  {\it (c)},  as case  $ \, n = 0 \, $.

\salto

\proclaim{Theorem 5.5}  Let  $ \, x = n \in \N \, $,  and let  $ \, f
\in \N_+ \, $,  $ \, k := \left[ {{f-n+1} \over {2}} \right] \, $.  Then
                                        \par
   \hskip-3pt  (a) no minor of order  $ (n+1) $  in  $ \, \Bfor
\setminus \Bfor(k) \, $  belongs to  $ \Rad \left( \Bfor \right)
\, $,
                                        \par
   \hskip-3pt  (b) \, no  $ f $--diagram  in  $ \, \Bfze \setminus
\Bfze \! \left( \big[ (f+1)/2 \big] \right) \, $  belongs to
$ \Rad \left( \Bfze \right) \, $,
                                        \par
   \hskip-3pt  (c) no Pfaffian of order  $ 2(n+1) $  in  $ \, \Bfsp
\setminus \Bfsp(k) \, $  belongs to  $ \Rad \left( \Bfsp \!\right)
\, $.
\endproclaim

\demo{Proof} {\it (a)} \, Let  $ \delta_{n+1} $  be a minor of order
$ (n+1) $  in  $ \, \Bfor(h) \setminus \Bfor (h+1) \, $,  \, with  $ \, h
< k \; $;  \, then  $ \delta_{n+1} $  is an alternating sum of diagrams
$ \d_i $  ($ i= 1, 2, \dots, (n+1)! \, $),  and at least one of these
diagrams   --- say $ \d_1 $  ---   has  $ h $  arcs, so it has at
least  $ n+1 $  vertical edges.  Let  $ v_1 $,  $ v_2 $,  $ \dots $,
$ v_{2(n+1)} $  be the moving vertices of  $ \delta_{n+1} \, $:  then
some of them, say  $ v_1 $,  $ v_2 $,  $ \dots $,  $ v_r $,  lay on the
top row, and the others, namely  $ v_{r+1} $,  $ v_{r+2} $,  $ \dots $,
$ v_{2(n+1)} $,  lay on the bottom row.  Here  $ \, 0 \leq r \leq 2(n+1)
\, $;  setting  $ \, s:= 2(n+1) - r \, $,  we can assume (by symmetry)
$ \, s \geq r \, $.  As  $ \, s+r = 2(n+1) \, $  is even,  $ s-r $
is even too, say  $ \, s-r = 2 \, t $,  $ \, t \in \N $.
                                                 \par
   We assumed that  $ \, s \geq r \, $,  but we reduce at once to the case
$ \, s=r \, $.  In fact, suppose  $ \, s \gneqq r \, $,  so  $ \, t \geq 1
\, $:  then in  $ \d_1 $  there are at least  $ t $  arcs on the bottom row
which pairwise join  $ 2 \, t $  vertices among  $ v_{r+1} $,  $ v_{r+2} $,
$ \dots $,  $ v_{2(n+1)} \, $.  Moreover, there are at least  $ r $
vertical edges joining  $ r $  vertices among  $ v_1 $,  $ v_2 $,
$ \dots $,  $ v_r $  with  $ r $  vertices among  $ v_{r+1} $,  $ v_{r+2} $,
$ \dots $,  $ v_{2(n+1)} $.  Since  $ \d_1 $  has exactly  $ \, f-2h \, $
vertical edges, we get that at least  $ \, f-2h-r \, $  of these vertical
edges do not involve  $ v_1 $,  $ v_2 $,  $ \dots $,  $ v_{2(n+1)} $,  hence
they are fixed in  $ \delta_{n+1} $,  i.e.~they also appear in all the other
summands  $ \d_i $  of  $ \delta_{n+1} \, $.  Now,  $ \, r+s = 2(n+1) \, $
and  $ \, s-r = 2t \, $  imply  $ \, r+t = n+1 \, $;  since  $ \, f-2h \geq
n+1 \, $,  we find  $ \, f-2h-r \geq n+1-r = t \, $;  therefore  $ \d_1 $
has at least  $ t $  vertical edges which also appear as well in all the
other summands  $ \d_i \, $.
                                             \par
   Let now  $ \h $  be a diagram selected as follows: pick  $ t $
vertices  $ w_1 $,  $ w_2 $,  $ \dots $,  $ w_t $  in the bottom row which in
$ \delta_{n+1} $  belong to fixed vertical edges, and pick  $ t $  vertices
$ v_{i_1} $,  $ v_{i_2} $,  $ \dots $,  $ v_{i_t} $  among  $ v_{r+1} $,
$ v_{r+2} $,  $ \dots $,  $ v_{2(n+1)} $   which are not joined with each
other in  $ \d_1 \, $.  Then let  $ \h $  be any  $ f $--diagram  which
has  $ t $  arcs in the top row joining each  $ w_i $  to a  $ v_{i_j}
\, $,  has  $ t $  arcs in the top row in the same positions than in
the bottom one, and whose remaining  $ f-2t $  vertices are joined by
straight-vertical edges (such an  $ \h $  is nothing but a suitable
product of  $ \h_{i,j} $'s).
%
% For instance, for  $ f=9 $,  $ t=3 $,
% $ \{w_1, \dots, w_t\} = \{5,6,9\} \, $,  $ \, \{v_{i_1}, \dots,
% v_{i_t}\} = \{1,3,7\} \, $  we have
% %
% %   PICTURE OF A SAMPLE OF A 9-DIAGRAM \h
% %
% \vskip5pt
%   $$  \raise27pt\hbox{\h \hskip23pt =}
% \hskip23pt  \epsfbox{eta-diagram.eps}  $$
% \vskip5pt
% %
% \indent
%
                                                         \par
   Now consider  $ \, \delta_{n+1} \cdot \h \, $.  By construction (see
Lemma 5.4 and its proof) we have that  $ \, \delta_{n+1}' := \delta_{n+1}
\cdot \h \, $  is a new minor (of order  $ n+1 $),  with again  $ \,
\delta_{n+1}' := \delta_{n+1} \cdot \h \in \Bfor(h) \, $.  But now
$ \delta_{n+1}' $  has  $ \, r'=r+t \, $  moving vertices in the top
row, and  $ \, s'=s-t \, $  moving vertices in the bottom row, thus  $ \,
r'=n+1=s' \, $.  Since  $ \Rad \left( \Bfor \right) $  is an ideal,
proving  $ \, \delta_{n+1}' \notin \Rad \left( \Bfor \right) \, $  will
also imply  $ \, \delta_{n+1} \notin \Rad \left( \Bfor \right) \, $,
as we wish to show.
                                                    \par
   Let us consider the  $ (n+1) $  moving vertices in the top row of the
minor  $ \delta_{n+1}' $.  As explained in  Remarks  4.7{\it(b)},  we
can split them into two disjoint subsets, that of  {\sl row vertices},
say  $ p $  in number, and that of  {\sl column vertices},  say  $ q $
in number (with  $ \, p+q = n+1 \, $).  In the bottom row of course we
have instead exactly  $ q $  row (moving) vertices and  $ p $  column
(moving) vertices.  Among the  $ f $--diagrams  in the expansion of
$ \delta_{n+1}' \, $,  we collect those attached to permutations in
$ \, S_p \times S_q \, $  ($ \, \subseteq S_{p+q} = S_{n+1} \, $),
i.e.~those in which the edges in the moving part are all vertical (which
act separately on the row and the column vertices in top row), and we denote
by  $ \varDelta $  their sum (with   --- alternating ---   signs).  Thus
we find that  $ \; \delta_{n+1}' = \varDelta + \varGamma \, $,  \; where
$ \varGamma $  is an algebraic sum of diagrams which all have the same
fixed part as  $ \varDelta $,  and moving part having at least one arc more
than  $ \varDelta \, $.  Namely, we have  $ \, \varDelta \in \Bfor \langle
h \rangle \, $  and  $ \, \varGamma \in \Bfor(h+1) \, $.  Therefore
 \vskip-6pt
  $$  \delta_{n+1}' \equiv \varDelta \, \mod \Bfor(h+1)  $$
 \vskip-3pt
\noindent
 so  $ \, \overline{\delta_{n+1}'} = \overline{\varDelta} \, $  as cosets
in  $ \Bfor[h] \, $.  Furthermore, if  $ \, \pi_1 \big( \delta_{n+1}' \big)
= 0 \, $   --- i.e.~$ \, \delta_{n+1}' \in \Rad \left( \Bfor \right) \, $
---   then clearly  $ \, \pi^*_1 \left(\, \overline{\varDelta} \,\right)
= \overline{0} \, $,  \, so it is enough to show that the latter cannot
occur.
                                                     \par
   By construction, the ``symmetric group part'' (i.e.~the one made of
vertical edges) of  $ \varDelta $  is just a product of some  $ \, \sigma
\, (\, \in S_f \,) $  times a product of antisymmetrizers  $ \, \hbox{\it
Alt}_p \cdot \hbox{\it Alt}_q \, $.  Since  $ \, \hbox{\it Alt}_p \cdot
\hbox{\it Alt}_q \, $  inside the group algebra  $ \Bbbk \big[ S_{f-2h}
\big] $  generates the two-sided ideal  $ \, \bigoplus_{\!\!\Sb
                  \! \mu \vdash (f-2h)  \\
                  \! \mu_1^t \geq p \, , \, \mu_2^t \geq q  \\
                \endSb} \hskip-5pt I_\mu \, $
(assuming  $ \, p \geq q \, $,  \, say), we can conclude that
$ \overline{\varDelta} $  generates the whole $ \Bfor $--bimodule
$ \; \bigoplus_{\Sb
                  \! \mu \vdash (f-2h)  \\
                  \! \mu_1^t \geq p \, , \, \mu_2^t \geq q  \\
                \endSb} \hskip-3pt \Bfor[h \, ; \mu] \; $
(cf.~\S 2).  Therefore, all of the
$ \Bfor $--bimodule  $ \; \bigoplus_{\!\!\Sb
                    \! \mu \vdash (f-2h)  \\
                    \! \mu_1^t \geq p \, , \, \mu_2^t \geq q  \\
                      \endSb} \hskip-3pt \Sfor[h \, ; \mu] \; $
is generated by  $ \, \pi_1^* \left(\, \overline{\varDelta} \,\right) \, $.
Then if  $ \, \pi_1^* \left(\, \overline{\varDelta} \,\right) = \overline{0}
\, $  we have also  $ \; \bigoplus_{\Sb
                    \!\!\! \mu \vdash (f-2h)  \\
                    \! \mu_1^t \geq p \, , \, \mu_2^t \geq q  \\
                      \endSb} \hskip-5pt \Sfor[h \, ; \mu] = 0 \, $;
\; as the latter is false   --- thanks to Proposition 3.2 ---   we
must have  $ \, \pi_1^* \left(\, \overline{\varDelta} \,\right) \neq
\overline{0} \, $,  q.e.d.
                                                \par
   {\it (b)} \, We can repeat the proof of the previous case, with the
obvious, wide simplifications (only the last third of that proof is
still necessary ---   and sufficient!).
                                                \par
   {\it (c)} \, The proof is similar to that of case  {\it (a)}.  Let
$ \varpi_{n+1} $  be a Pfaffian of order  $ 2(n+1) $  in  $ \, \Bfsp(h)
\setminus \Bfsp(h+1) \, $,  for some  $ \, h < k \; $;  \, then
$ \varpi_{n+1} $  is a sum of diagrams  $ \d_i $  ($ i= 1, 2, \dots,
(n+1)! \, $),  and at least one of them   --- say  $ \d_1 $  ---
has  $ h $  arcs, hence it has at least  $ n+1 $  vertical edges.  Let
$ v_1 $,  $ v_2 $,  $ \dots $,  $ v_{2(n+1)} $  be the moving vertices
of  $ \varpi_{n+1} \, $:  then some of them, say $ v_1 $,  $ v_2 $,
$ \dots $,  $ v_r $,  lay on the top row, and the others, namely
$ v_{r+1} $,  $ v_{r+2} $,  $ \dots $,  $ v_{2(n+1)} $,  lay on the
bottom row (with  $ \, 0 \leq r \leq 2(n+1) \, $,  $ \, s:= 2(n+1)
- r \, $).  Exactly as in  {\it (a)},  we can reduce to the case
$ \, s=r \, (=n+1) \, $.
                                                    \par
   Now the Pfaffian  $ \varpi $  has  $ \, n+1 \, $  moving vertices up
and  $ \, n+1 \, $  down.  Hence among the diagrams in the sum expressing
$ \varpi $  there are some whose moving edges are all vertical: namely,
those corresponding to the terms in (4.3) with  $ \, \{ h_i \mid i=1, 2,
\dots, n+1 \} \subseteq \{1,2, \dots, f \,\} \, $  (the top row) and  $ \,
\{k_1, k_2, \dots, k_{n+1} \} \subseteq \{f+1,f+2, \dots, 2f \,\} \, $  (the
bottom row).  But  $ h_1 $,  $ h_2 $,  $ \dots $,  $ h_{n+1} $  are fixed
by the condition  $ \, h_1 < h_2 < \cdots < h_{n+1} \, $,  whilst there is
no condition on the ordering of  $ k_1 $,  $ k_2 $,  $ \dots $,  $ k_{n+1}
\, $.  Thus all diagrams of the previous type are obtained by fixing the
sets  $ \, \{h_1, h_2, \dots, h_{n+1} \} \, $  and  $ \, \{k_1, k_2, \dots,
k_{n+1} \} \, $  and joining  $ h_i $  to  $ k_{\sigma(i)} $  for all  $ \,
i=1, 2, \dots, n+1 $,  for all  $ \, \sigma \in S_{n+1} \, $.  We denote
by  $ \varPi $  the sum of these diagrams, and we note that the
$ S_{n+1} $--action  on  $ \, \{h_1, h_2, \dots, h_{n+1} \} \, $  or 
$ \, \{k_1, k_2, \dots, k_{n+1} \} \, $  turns  $ \varpi $  into itself.
                                                     \par
   So far we found that  $ \; \varpi_{n+1} = \varPi + \varGamma \, $,  \;
where  $ \varGamma $  is a sum of diagrams which all have the same fixed
part as  $ \varPi $  and moving part having at least one arc more than
$ \varPi \, $.  That is, we have  $ \, \varPi \in \Bfsp \langle h \rangle
\, $  and  $ \, \varGamma \in \Bfsp(h+1) \, $.  Thus
 \vskip-8pt
  $$  \varpi_{n+1} \equiv \varPi \, \mod \Bfsp(h+1)  $$
 \vskip-4pt
\noindent
 and so  $ \; \overline{\varpi_{n+1}} = \overline{\varPi} \in \Bfsp[h]
\, $.  Moreover, here again  $ \, \varpi'_{n+1} \in \Rad \left( \Bfsp
\right) \, $  would imply  $ \, \pi_1^* \left(\, \overline{\varPi} \,
\right) = \overline{0} \, $  as well, thus we have to show that the
latter does not occur.
                                                     \par
   By construction, the ``symmetric group part'' of  $ \varPi $  is
just the product of some  $ \, \sigma \, (\, \in S_f \,) $  times a
symmetrizer  $ \, \hbox{\it Sym}_{n+1} \, $,  whence it follows that
$ \overline{\varPi} $  generates the whole  $ \Bfsp $--bimodule
$ \; \bigoplus_{\Sb
                  \! \mu \vdash (f-2h)  \\
                  \! \mu_1 \geq (n+1)  \\
                \endSb} \hskip-2pt \Bfsp[h \, ; \mu] \; $.
Therefore,  $ \pi_1^* \left(\, \overline{\varPi} \,\right) $  in turn
generates the  $ \Bfsp $--bimodule  $ \; \bigoplus_{\Sb
                         \! \mu \vdash (f-2t)  \\
                         \! \mu_1 \geq (n+1)  \\
                        \endSb} \hskip-2pt \Sfsp[h \, ; \mu] \, $,
\; hence one has that  $ \, \pi_1 \big( \varPi \big) = 0 \, $  would imply
also that  $ \; \bigoplus_{\Sb
                 \! \mu \vdash (f-2t)  \\
                 \! \mu_1 \geq (n+1)  \\
                \endSb} \hskip-2pt \Sfsp[h \, ; \mu] = 0 \; $;
\; the latter is false, by Proposition 3.2, so  $ \, \pi_1^*
\left(\, \overline{\varPi} \,\right) \neq 0 \, $.   $ \square $
\enddemo

\vskip3pt

   Theorems 5.3 and 5.5 have many consequences for the radicals of the
various algebras  $ \Bfx(h) $  and  $ \Bfx[h] \, $,  \, which we collect
in the following two statements.  Here again, we wish to point out that,
altogether, these corollaries give a necessary and sufficient condition
for the coset of a minor   --- of order  $ (n+1) $  ---   respectively
of an  $ f $--diagram,  respectively of a Pfaffian   ---   of order
$ 2(n+1) $  ---   to belong to the radical of the suitable quotient
algebra.
%
% Note also that  Corollary 5.6{\it (b)\/}  extends to the
% case of {\sl any}  $ f $  what is proved in Proposition 3.2 for
% the case of  {\sl even}  $ f $  only.
%

\salto

\proclaim{Corollary 5.6}  Let  $ \, n \in \N_+ \, $,  $ \, f \in \N_+ \, $,
$ \, k := \big[ {{f-n+1} \over {2}} \big] \, $,  \, and  $ \, h \geq k
\; $.  \, Then
                                                 \par
   (a) \, every minor of order  $ (n+1) $  in  $ \Bfor\!(h) $  belongs
to  $ \Rad \left( \Bfor\!(h) \right) \, $,  \, and its coset either in
$ \, \Bfor \! \Big/ \Bfor\!(h+1) \, $  or in  $ \, \Bfor[h] := \Bfor\!(h)
\! \Big/ \Bfor\!(h+1) \, $  belongs to  $ \, \Rad \left( \Bfor \! \Big/
\Bfor\! (h+1) \right) $  or  $ \, \Rad \left( \Bfor[h] \right) \, $, 
\, respectively;
                                                 \par
   (b) \, every  $ f $--diagram  in  $ \Bfze \! \left( \big[ (f+1)/2 \big]
\right) \, $  belongs to  $ \Rad \left( \Bfze \! \left( \big[ (f+1)/2 \big]
\right) \right) \, $,  \, hence  $ \; \Bfze \! \left( \big[ (f+1)/2 \big]
\right) = \Rad \left( \Bfze \! \left( \big[ (f+1)/2 \big] \right) \right)
\, $,  \; so that  $ \, {\Cal S}^{(0)}_f \! \left( \big[ (f+1)/2 \big]
\right) = 0 \; $;
                                                 \par
   (c) \, every Pfaffian of order  $ \, 2(n+1) $  in  $ \Bfsp\!(h) $
belongs to  $ \Rad \left( \Bfsp\!(h) \!\right) \, $,  \, and its coset
either in  $ \, \Bfsp \! \Big/ \Bfsp\!(h+1) \, $  or in  $ \, \Bfsp[h]
:= \Bfsp\!(h) \! \Big/ \Bfsp\!(h+1) \, $  belongs to  $ \, \Rad \left(
\Bfsp \! \Big/ \Bfsp\! (h+1) \right) $  or  $ \, \Rad \left( \Bfsp[h]
\right) \, $,  \, respectively.
\endproclaim
 \eject

\demo{Proof} Everything follows directly from Theorem 5.2 once we remind
that
  $$  \Rad\,(R) \; {\textstyle \bigcap} \; I \,\; = \; \Rad\,(I\,)  \quad ,
\qquad  \big( \Rad\,(R) + I \,\big) \Big/ I \,\; = \; \Rad \Big( R \big/
I \Big)  $$
for every ideal  $ I $  in an Artinian ring  $ R $   --- as observed in
[HW1] ---   and we apply this fact to the cases  $ \, R = \Bfx \, $  or
$ \, R = \Bfx(h) \, $  and  $ \; I = \Bfx(h+1) \;\; $.   \qed
\enddemo

\salto

\proclaim{Corollary 5.7}  Let  $ \, n \in \N_+ \, $,  $ \, f \in \N_+ \, $,
$ \, k := \big[ {{f-n+1} \over {2}} \big] \, $,  \, and  $ \, h < k \; $.
                                                 \par
   (a) \, let  $ \delta_{n+1} $  be a minor of order  $ (n+1) $  in
$ \Bfor(h) \, $.  Then the coset of  $ \, \delta_{n+1} $  either in
$ \, \Bfor \! \Big/ \Bfor\!(h+1) \, $  or in  $ \, \Bfor[h] := \Bfor\!(h)
\! \Big/ \Bfor\!(h+1) \, $  does  {\sl not}  belong to the radical  $ \,
\Rad \left( \Bfor \! \Big/ \Bfor\! (h+1) \right) $  nor to  $ \, \Rad
\left( \Bfor[h] \right) $  respectively;
                                                 \par
   (b) \, let  $ \d $  be an  $ f $--diagram  in  $ \, D_{f,h} \setminus
D_{f,[(f+1)/2]} \; $.  Then the coset of  $ \, \d $  either in  $ \, \Bfze
\! \Big/ \Bfze\!\big([f/2]+1\big) \, $  or in  $ \, \Bfze[h] := \Bfze\!(h)
\! \Big/ \Bfze\!(h+1) \, $  does  {\sl not}  belong to the radical  $ \,
\Rad \left( \Bfze \! \Big/ \Bfze\!\big([f/2]+1\big) \right) $  nor to
$ \, \Rad \left( \Bfze[h] \right) $  respectively;
                                                 \par
   (c) \, let  $ \varpi_{n+1} $  be a Pfaffian of order  $ \, 2(n+1) $  in
$ \, \Bfsp(h) \, $.  Then the coset of  $ \, \varpi_{n+1} $  either in  $ \,
\Bfsp \! \Big/ \Bfsp\!(h+1) \, $  or in  $ \, \Bfsp[h] := \Bfsp\!(h) \!
\Big/ \Bfsp\!(h+1) \, $  does  {\sl not}  belong to the radical  $ \,
\Rad \left( \Bfsp \! \Big/ \Bfsp\! (h+1) \right) $  nor to  $ \, \Rad
\left( \Bfsp[h] \right) $  respectively.
\endproclaim

\demo{Proof} The claim follows easily from the analysis carried out to
prove Theorem 5.5.   \qed
\enddemo

\salto

   {\bf 5.8 An overall conjecture about the radical.} \  Theorem 5.3
provides the following global information about the radical of Brauer
algebras:

\vskip5pt

   \hskip-5pt  {\it   --- a) \hskip11pt  $ \Rad \left( \Bfor \right) \,
\supseteq \, R^{(n)}_f \, := \; \text{\sl  $ \Bbbk $--span}\bigg( \text{\sl
Min}^{\,(n)}_{f;\,n+1} \, \bigcap \, \Bfor \! \Big( \big[ {(f \! -
\! n \! + \! 1) \big/ 2} \big] \Big) \bigg) $}

\vskip3pt

   \hskip-5pt  {\it   --- b) \hskip11pt  $ \Rad \left( \Bfze \right)
\, \supseteq \, R^{(0)}_f \, := \; \text{\sl  $ \Bbbk $--span}\Big(
D_{f, \, [{(f+1) / 2}]} \Big) \, = \, \Bfze \! \left( \big[ (f+1)
\big/ 2 \big] \right) $}

\vskip3pt

   \hskip-5pt  {\it   --- c) \hskip11pt  $ \Rad \left( \Bfsp \right) \,
\supseteq \, R^{(-2n)}_f \, := \; \text{\sl  $ \Bbbk $--span}\bigg(
\text{\sl Pf}^{\;(-2n)}_{f;\,n+1} \, \bigcap \, \Bfsp \! \Big( \big[
{(f \! - \! n \! + \! 1) \big/ 2} \big] \Big) \bigg) $}

\vskip7pt

   In principle,  $ \Rad \left( \Bfx \right) $  might be greater than the
space  $ R^{(x)}_f $  considered above.  Nevertheless, at least in some
cases we are able to leave out of the radical some elements which,  {\it
a priori},  might belong to it.  We shall now briefly sketch what we mean.
                                                 \par
   Let us consider for instance the case of  $ \Bfor $,  with  $ \, n \in
\N_+ \, $.  Pick  $ \, \eta \in \! \Rad \left( \Bfor \right) \, $.  Resuming
notations of \S 5.1, we have  $ \, \eta \in \! \Rad \left( \Bfor \right)
\subseteq \text{\it Ker}\big(\pi_V\big) = \text{\sl  $ \Bbbk $--span} \Big(
\text{\sl Min}^{\,(n)}_{f;\,n+1}\Big) \, $.  Now, fix a subset  $ B $  of
$ \, \text{\sl Min}^{\,(n)}_{f;\,n+1} \, $  which is a  $ \Bbbk $--basis
of  $ \text{\it Ker}\big(\pi_V\big) \, $,  \, and expand  $ \eta $  as
a  $ \Bbbk $--linear  combination of elements of  $ B \, $.  Since  $ \;
\text{\sl Min}^{\,(n)}_{f;\,n+1} \, \bigcap \, \Bfor \! \Big( \big[ {(f \!
- \! n \! + \! 1) \big/ 2} \big] \Big) \; $  is contained in  $ \Rad \left(
\Bfor \right) $,  we can even reduce to the case where in such an expansion
of  $ \eta $  there occur (with non-zero coefficient)  {\sl only elements
of}  $ \, B \setminus \Bfor \! \Big( \big[ {(f \! - \! n \! + \! 1) \big/ 2}
\big] \Big) \, $.
%
%                                                  \par
%
 \eject
   {\sl In particular, let us assume, for instance, that  $ \, \eta = c_1
\, \delta_1 + c_2 \, \delta_2 \, $},  \, where we have  $ \, \delta_1 ,
\delta_2 \in B \setminus \Bfor \! \Big( \big[ {(f \! - \! n \! + \! 1)
\big/ 2} \big] \Big) \, $,  $ \, \delta_1 \not= \delta_2 \, $  and  $ \,
c_1 , c_2 \in \Bbbk \, $;  {\sl in addition, we assume that  $ \delta_2 $
has two distinct moving vertices on a single row, say  $ i^+ $  and
$ j^+ $,  which are  {\it not}  both moving for  $ \delta_1 \, $,
\, nor both always on (fixed or ``moving'') vertical edges.}
                                                 \par
   In this situation, we consider  $ \, \h_{i,j} \cdot \eta \, $,  \, which
clearly also belongs to  $ \Rad \left( \Bfor \right) \, $.  Then Lemma 5.4
gives  $ \; \h_{i,j} \cdot \delta_2 = 0 \; $  and  $ \; \h_{i,j} \cdot
\delta_1 = n^e \, \delta'_1 \; $  with  $ \, e \in \N \, $  and  $ \,
\delta'_1 \in \text{\sl Min}^{\,(n)}_{f;\,n+1} \, $,  \; so

 \vskip-5pt
  $$  \Rad \left( \Bfor \right)  \, \ni \;  \h_{i,j} \cdot \eta \, = \,
c_1 \, \h_{i,j} \cdot \delta_1 + c_2 \, \h_{i,j} \cdot \delta_2 \, = \,
c_1 \, n^e \delta'_1  \quad .  $$
 \vskip-5pt
\noindent
 But then it follows that  $ \, \delta'_1 \in \Rad \left( \Bfor \right)
\, $,  \, and also  $ \, \delta'_1 \not\in \Bfor \! \Big( \big[ {(f \!
- \! n \! + \! 1) \big/ 2} \big] \Big) \, $,  \, by construction, which
by  Theorem 5.5{\it (a)\/}  is impossible.
                                                 \par
   Similar results can be obtained (via  Theorem 5.5{\it (c)\/})  in the
case of  $ \Bfsp $  ($ n \! \in \! \N_+ $)  too.
                                                 \par
   On the other hand, in some cases the inclusion of  $ R_f^{(z)} $  in
$ \Rad \left( {\Cal B}_f^{(z)} \right) $  is strict; for instance, this
is the case for  $ \, z = 0 \, $  and  $ \, f \geq 7 \, $,  \, by
definitions and by the results in [Ru] and [RS], for which  $ \, {\Cal B}_f^{(0)} $  is semisimple iff  $ \, f \in \{1,3,5\} \; $.

\vskip3pt

   However, we can also remark that (cf.~the proof of Corollary 5.6)

\vskip3pt

   \centerline{ $ R^{(n)}_f  \; \subseteq \;
\Rad \left( \Bfor \right) \bigcap \, \Bfor \! \Big( \big[ {(f \!
- n + \! 1) \big/ 2} \big] \Big)  \, = \,  \Rad \, \bigg( \Bfor \!
\Big( \big[ {(f \! - n + \! 1) \big/ 2} \big] \Big) \! \bigg) $ }

\vskip2pt

   \centerline{ $ R^{(0)}_f  \; \subseteq \;
\Rad \left( {\Cal B}_f^{(0)} \right) \bigcap \, {\Cal B}_f^{(0)}
\! \Big( \big[ {(f \! + \! 1) \big/ 2} \big] \Big)  \, = \,
\Rad \, \bigg( {\Cal B}_f^{(0)} \! \Big( \big[ {(f \! + \! 1)
\big/ 2} \big] \Big) \! \bigg) $ }

\vskip2pt

   \centerline{ $ R^{(-2n)}_f  \; \subseteq \;
\Rad \left( \Bfsp \right) \bigcap \, \Bfor \! \Big( \big[ {(f \!
- n + \! 1) \big/ 2} \big] \Big)  \, = \,  \Rad \, \bigg( \Bfsp \!
\Big( \big[ {(f \! - n + \! 1) \big/ 2} \big] \Big) \! \bigg) $ }

\vskip3pt

\noindent
 Even more, the mid-line inclusion is indeed an identity, namely

\salto
\vskip3pt

\proclaim{Proposition 5.9}  Let  $ \, f \in \N_+ \, $.  Then
$ \;\; \displaystyle{ \Rad \, \Big( {\Cal B}_f^{(0)} \! \Big( \big[ {(f \! + \! 1)
\big/ 2} \big] \Big) \Big)  \, = \,  R^{(0)}_f } \; $.
%
%   $$  \Rad \, \Big( {\Cal B}_f^{(0)} \! \Big( \big[ {(f \! + \! 1)
% \big/ 2} \big] \Big) \Big)  \; = \;  R^{(0)}_f  $$
%
\endproclaim

\demo{Proof}  If  $ f $  is odd, then  $ \, {\Cal B}_f^{(0)} \!
\Big( \big[ {(f \! + \! 1) \big/ 2} \big] \Big) = 0 \, $  and  $ \,
R^{(0)}_f = 0 \, $  by definition.  If  $ f $  is even, then  $ \,
R_f^{(0)} := {\Cal B}_f^{(0)} \big( {f \big/ 2} \big) = {\Cal B}_f^{(0)}
\! \Big( \big[ {(f \! + \! 1) \big/ 2} \big] \Big) \, $,  \, and the
claim follows from Corollary 3.3.   \qed
\enddemo

\saltino

   Thus, inspired by these experimental evidences, we are lead to
formulate the following

\ssalto

\proclaim{Conjecture 5.10}  Let  $ \, f \in \N_+ \, $,  $ \, n \in \N_+
\, $.  Then
 \vskip-17pt
  $$  \Rad \left( \Bfor \! \Big( \big[ {(f \! - \! n \! + \! 1) \big/ 2}
\big] \Big) \right) \; = \; R^{(n)}_f  \quad ,  \qquad
  \Rad \left( \Bfsp \! \Big( \big[ {(f \! - \! n \! + \! 1) \big/ 2}
\big] \Big) \right) \; = \; R^{(-2n)}_f  $$
\endproclaim

\vskip5pt

   {\bf 5.11 Inheriting the radical.} \  In [HW2], \S 3.2, an algorithm is
described for constructing a part of  $ \Rad \left( \Bfx \right) $,  \,
called ``the Inherited Piece of the Radical'', out of  $ \Rad \left(
{\Cal B}^{(x)}_{f-2} \right) \, $.
                                             \par
  The construction is the following.  Take an  $ (f\!-\!2) $--diagram  $ \,
\d \in D_{f-2} \, $,  and let  $ (i,j) $  and  $ (h,k) $  be two pairs of
numbers such that  $ \, 1 \leq i < j \leq f \, $  and  $ \, 1 \leq h < k
\leq f \, $.  Define an  $ f $--diagram  $ \, \d^{i,j}_{h,k} \in D_f \, $
to be the diagram obtained from  $ \d $  by inserting a new arc \  $ \,
i^+ \epsfbox{r--s_barup.eps} \, j^+ \, $  \ in the top row and a new arc
\  $ \, h^- \epsfbox{r--s_bardown.eps} \, k^- \, $  \ in the bottom row.
Then let  $ \; {Ex}^{i,j}_{h,k} : \, {\Cal B}^{(x)}_{f-2} \longrightarrow
\Bfx \, $  be the unique linear embedding defined by  $ \, {Ex}^{i,j}_{h,k}
(\d) := \d^{i,j}_{h,k} \, $  for all  $ \, \d \in D_{f-2} \, $.  Given a
subspace  $ \, {\Cal I} \, $  of  $ {\Cal B}^{(x)}_{f-2} \, $,  \, define
$ \, {\Cal I}^{\,(1)} := \sum_{(i,j),(h,k)} {Ex}^{i,j}_{h,k} \big(
{\Cal I} \,\big) \, $,  \, the  $ \Bbbk $--span  of all the
$ {Ex}^{i,j}_{h,k} \big( {\Cal I} \,\big) $'s.  Then Theorem
3.2.9 in [HW2] claims that

\vskip5pt

  \centerline{ \it  $ {\left( \Rad \left( {\Cal B}^{(x)}_{f-2} \right)
\! \right)}^{\!(1)} $  is a two-sided ideal of  $ \, \Bfx $,  and it
is contained in  $ \Rad \left( \Bfx \right) \, $. }

\vskip5pt

   Actually, this also follows (see [KX], Lemma 3.1) because  $ \Bfx $
is a cellular algebra.

\vskip3pt

  Now the remark is the following.  Let  $ \eta $  be an element of
$ \Rad \left( {\Cal B}_{f-2}^{(x)} \right) $  of the type given in
Theorem 5.3 (i.e.~a minor or a Pfaffian).  By the very definitions,
$ \, {Ex}^{i,j}_{h,k} (x) \, $  is again an element of the same type,
so Theorem 5.3 applies to give  $ \, {Ex}^{i,j}_{h,k}(\eta) \in \Rad
\left( \Bfx \right) \, $.  Therefore, if  $ {\Cal I} $  is the span
of the previous elements in  $ {\Cal B}_{f-2}^{(x)} \, $,  then we see
directly that  $ \, {\Cal I}^{\,(1)} \subseteq \Rad \left( \Bfx \right)
\, $,  \, which is a nice way to verify the ``inheritance phenomenon''
of [HW2].
                                                  \par
  Note also that this remark is fully consistent with Proposition 5.9
and Conjecture 5.10.

\salto
 \vskip7pt

   {\bf 5.12 The radical of the  $ \Bfx $--modules  $ H^\mu_{f,k} \, $.}
\  The information about the radical of a Brauer algebra  $ \Bfx $
given in Theorem 5.3 yield also information about the radical of the
$ \Bfx $--modules  $ H^\mu_{f,k} \, $.  Namely, the first, immediate
result is (with notation of \S 5.8)

\salto
 \vskip5pt

\proclaim{Theorem 5.13} Let  $ \, n \in \N_+ \, $,  $ \, f \in \N_+ \, $,
$ \, k \in \big\{ 0, 1, \dots, [f/2] \big\} \, $,  and  $ \, \mu \vdash
( f - 2k \,) \, $.  Then
 \vskip-15pt
  $$  R^{(n)}_f.\,H^\mu_{f,k} \, \subseteq \, \Rad\, \big( H^\mu_{f,k}
\big)  \;\; ,  \quad  R^{(0)}_f.\,H^\mu_{f,k} \, \subseteq \, \Rad\,
\big( H^\mu_{f,k} \big)  \;\; ,  \quad  R^{(-2n)}_f.\,H^\mu_{f,k}
\, \subseteq \, \Rad\, \big( H^\mu_{f,k} \big)  $$
 \vskip-5pt
\noindent
 where in each case  $ H^\mu_{f,k} $  stands for the suitable module
for the algebra  $ \Bfor $,  $ \Bfze $  or  $ \Bfsp $.
\endproclaim

\demo{Proof}  The claim follows immediately from the standard inclusion
$ \; \Rad\,(A).M \subseteq \Rad\,(M) \, $,  \; which holds for every ring
$ A $  and every  $ A $--module  $ M $.   \qed
\enddemo

\salto

   In addition, Conjecture 5.10 about the radical of the Brauer
algebras also involve a similar conjecture about the radicals of
the  $ H^\mu_{f,k} $'s,  namely that the inclusions in Theorem 5.13
actually be identities.  We prove this conjecture in some special cases,
see \S 6 below; also, the anlogous claim for  $ \Bfze $  is easily seen
to be true.  The complete statement is

\salto
\vskip3pt

\proclaim{Conjecture 5.14}  If  $ \, n \! \in \! \N_+ \, $,  $ \, f \! \in
\! \N_+ \, $,  $ \, k \in \! \big\{ [(f \! - \! n \! + \! 1)/2], \dots, [f/2] \big\} \, $,
$ \, \mu \vdash \! ( f \! - \! 2k ) \, $,  then
 \vskip-3pt
  $$  \Rad \, \big( H^\mu_{f,k} \big) \, = \, R^{(n)}_f.\,H^\mu_{f,k} \;\; ,
\hskip51pt
      \Rad \, \big( H^\mu_{f,k} \big) \, = \, R^{(-2n)}_f.\,H^\mu_{f,k}  $$
 \vskip-2pt
\noindent
 where in each case  $ H^\mu_{f,k} $  stands for the suitable module
for the algebra  $ \Bfor $  or  $ \Bfsp $.
\endproclaim

\demo{Proof}  If Conjecture 5.10 holds true, then the claim follows
immediately from the standard identity  $ \; \Rad\,(M) = \Rad\,(A).M
\, $,  \; which holds for every ring  $ A $  and every  $ A $--module
$ M $  such that  $ \, M \big/ \Rad\,(M) \, $  is simple, which is the
case for  $ \, M = H^\mu_{f,k} \; $,  \; by Proposition 2.7.   \qed
\enddemo

\salto

   {\sl  $ \underline{\text{Remark}} $:} \, in \S 6 below we shall see
some cases in which Conjecture 5.10 and Conjecture 5.14 do hold true,
namely when  $ \, n = 1 \, $  and  $ f $  is even.  For the sake of
completeness, we report also that [Ga], Corollary 4.6, gives other
results about  $ \Rad \left( \Bfx \right) $  and   $ \Rad \big(
H^\mu_{f,k} \big) \, $.

\ssalto

   {\bf 5.15 The case of positive characteristic.} \  Let  $ \, p :=
\text{\sl Char}\,(\Bbbk) \, $.  All our results about  $ \Rad \left(
\Bfx \right) $  in this section are based on the assumption  $ \, p
= 0 \, $.  We shall now discuss to what extent these results might
hold for  $ \, p > 0 \, $  as well.

\vskip2pt

   First of all, the results of \S 2 and \S 3 about the  $ \Bfx $--modules
and the semisimple quotient  $ \Sfx $  of  $ \Bfx $  only require  $ p > f
\; $.  On the other hand, the results of \S 4
%
% concerning the relations of
% Brauer algebras with centralizer algebras
%
 assume condition (4.1) to be
satisfied.  The latter does hold, in particular, whenever  $ \, p > f \; $;
\, therefore, under the latter, stronger assumption all results of \S\S 2,
3 and 4 are available.
                                                      \par
   On the other hand, from \S 5.1 on we assumed  $ \, p = 0 \, $.  This ensures
that the group  $ O(V) \, $,  resp.~$ Sp(W) \, $,  is linearly reductive: so
$ V^{\otimes f} $,  resp.~$ W^{\otimes f} $,  is a semisimple module for
$ O(V) \, $,  resp.~for  $ Sp(W) \, $.  Then, by general theory,  {\sl the centralizer algebra  $ \endor \, $,  resp.~$ \endsp \, $,  \, is semisimple
too}.  This last fact is the basis to obtain Proposition 5.2; the precise block decomposition) of this semisimple algebra then can be recovered from the results
of \S\S 2--4   --- see the proof of Proposition 5.2, where we summarized all this quoting [Wz].  Then Theorem 5.3, our first main result, follows as a direct consequence.

\vskip4pt

   We must point out a key fact.  Let  $ \, U \in \big\{ V, W \big\} \, $,  \,
let  $ \, G(V) := O(V) \, $,  $ \, G(W) := Sp(W) \, $.  For  $ \, E \subseteq {End}_{G(U)} \! \left( U^{\otimes f} \right) \, $  let  $ \, \big\langle E \big\rangle \, $  be the subalgebra  of  $ {End}_\Bbbk \! \left( U^{\otimes f} \right) $  generated by  $ U \, $.  Then

\ssalto

\proclaim{Lemma 5.16} Assume  $ \, p > f \, $.  Then the following are equivalent:
 \vskip2pt
   (a) \,  $ V^{\otimes f} $,  resp.~$ W^{\otimes f} $,  is a semisimple module for  $ O(V) \, $,  resp.~for  $ Sp(W) \, $;
                                                      \par
   (b) \,  the  $ \Bbbk $--algebra  $ \big\langle O(V) \big\rangle \, $,
resp.~$ \big\langle Sp(W) \big\rangle \, $,  is semisimple;
                                                      \par
   (c) \,  the  $ \Bbbk $--algebra  $ \pi_V\Big(\Bfor\Big) \, $,  \,
resp.~$ \pi_W\Big(\Bfsp\Big) \, $,  \, is semisimple;
                                                      \par
   (d) \,  $ V^{\otimes f} $,  resp.~$ W^{\otimes f} $,  is a semisimple module for  $ \pi_V\Big(\Bfor\Big) \, $,  resp.~for  $ \pi_W\Big(\Bfsp\Big) \, $.
\endproclaim

\demo{Proof} By Theorem 4.4,  $ \big\langle O(V) \big\rangle $  and  $ \, \pi_V \Big(\Bfor\Big) $  are the centralizer of each other, for their action on
$ V^{\otimes f} $,  and similarly for  $ \big\langle Sp(W) \big\rangle $  and
$ \, \pi_W\Big(\Bfsp\Big) $  acting on  $ W^{\otimes f} \, $.  Therefore, 
%   
% the Double Centralizer Theorem   
%   
 Schur duality    
tells us that  {\it (b)\/}  and  {\it (c)\/}  are equivalent.
                                                      \par
   On the other hand, the equivalences  \, {\it (a)\/}  $ \Longleftrightarrow $
{\it (b)\/} \,  and  \, {\it (c)\/}  $ \Longleftrightarrow $  {\it (d)\/} \,
are obvious.   \qed
\enddemo

\vskip3pt

   To sum up, the above analysis proves the following result:

\vskip7pt

\proclaim{Theorem 5.17}  Let  $ \, p := \text{\sl Char}\,(\Bbbk) > f \, $.
If any one of the (equivalent) conditions in Lemma 5.16 is satisfied, then
Proposition 5.2 and Theorem 5.3 still hold true.   \qed
\endproclaim

\vskip3pt

   On the other hand, the proof of Theorem 5.5 (our second main result) only
exploits combinatorial techniques and some results of \S 3.  Thus, its proof
still is valid if  $ \, p > f \, $;  \, so

\salto

\proclaim{Theorem 5.18}  Let  $ \, p := \text{\sl Char}\,(\Bbbk) > f \, $.
Then Theorem 5.5 still holds true.   \qed
\endproclaim

\vskip3pt

   Finally, we extend Conjectures 5.10 and 5.14 to the case of
$ \, p := \text{\sl Char}\,(\Bbbk) > f \, $,  \, too.
%
% \Salto
%
 \eject

\centerline{ \bf  \S 6  Applications: the Temperley-Lieb algebra and
pointed chord diagrams }

\ssalto

   {\bf 6.1  Temperley-Lieb algebra and pointed chord diagrams.} \  Let
$ \, f \in \N_+ \, $  be  {\sl even}.  In this section we study the cases
of  $ \Bfuno $  and  $ H^{(0)}_{f,f/2} \, $:  \, in particular, we compute
$ \Rad \left( \Bfuno\big(f/2\big) \right) $  and  $ \Rad \left( H^{(0)}_{f,f/2}
\right) \; $,  \, proving that Conjecture 5.10 and Conjecture 5.14 (respectively)
do hold true for them.  To this end, we shall not make use of Theorem 5.3; in
particular, we need no special assumptions on the ground field.

\vskip3pt

   Indeed,  {\bf in this section we assume that  $ \Bbbk $  is any field.}

\vskip5pt

   First, a terminological remark.  The  {\sl unital\/}  subalgebra of
$ \Bfuno $  generated by  $ \Bfuno\big(f/2\big) $  is usually called
{\sl Temperley-Lieb algebra\/}  (possibly defined in other ways, usually
involving some parameter too: see e.g.~[DN], and references therein), call
it  $ \Cal{TL}_f \, $.  The restriction functor yields an equivalence between
the category of all  $ \Cal{TL}_f \, $--modules and  the category of all 
$ \Bfuno\big(f/2\big) $--modules,  so studying the latter we are studying
the former too.
                                             \par
   Second, the set  $ J_{f,f/2} $  of all  $ \big( f, f/2 \big) $--junctions
(a  $ \Bbbk $--basis  of  $ H^{(0)}_{f,f/2} \, $)  can be represented by  {\sl
pointed chord diagrams},  as follows.  Given  $ \, j \in J_{f,f/2} \, $,  \,
let us lay the  $ f $  vertices of  $ j $  on the interior of a (horizontal)
segment, and draw the arcs of  $ j $  as arcs above this segment.  Now close
up the segment into a circle (like winding up the segment around a circle of
same length), gluing together the vertices of the segment and sealing them
with a special, marking dot.  Then the arcs of  $ j $  are turned into  {\sl
chords\/}  of the circle.  This (up to details) sets a bijection from
$ J_{f,f/2} $  to the set of pointed chord diagrams on the circle with
$ f $  chords.

\vskip3pt

   We'd better point out that ours are  {\sl pointed\/}  chord
diagrams.  Indeed, usually chord diagrams are considered up to
rotations, which is not the present case   ---   roughly, the marked
point forbids rotations.  Our previous construction shows that
one can identify  $ H^{(0)}_{f,f/2} $  with the  $ \Bbbk $--span
of the set of all pointed chord diagrams on the circle with  $ f $
chords: so this  $ \Bbbk $--span  is a module for  $ \Bfuno \big( f/2
\big) $  (hence for  $ \Bfuno $  as well)  or the Temperley-Lieb algebra.
It is customary to drop the ``pointed'' datum, considering chord diagrams
on the circle up to rotations: this amounts to look at the vertices of a
junction up to cyclical permutations.  In this case, the  $ \Bbbk $--span
of the set of all  {\sl non-pointed\/}  chord diagrams on the circle with
$ f $  chords (or of the set of all  $ \big(f,f/2\big) $--junctions  up
to cyclical permutations of their vertices) bears a natural structure of
module for the quotient algebra of  $ \Cal{TL}_f $   --- or of  $ \Bfuno
\big( f/2 \big) $,  or of  $ \Bfuno $  too ---   given by taking
$ f $--diagrams  up to simultaneous cyclic equivalence of their
top and bottom vertices.  This algebra is sometimes called  {\sl
Temperley-Lieb  {\rm (or, respectively,  {\sl Brauer\/})}  {\it
loop\/}   --- or  {\it affine}  ---   algebra}.  Indeed, this is
the most common framework where chord diagrams, and Temperley-Lieb
or Brauer algebras acting on them, do appear in literature: see, for
instance, [DN], or [Jo], and references therein.  The results we find
below for modules built upon pointed chord diagrams and for Temperley-Lieb
or Brauer algebras can be easily adapted to the non-pointed and the
loop/affine case as well.

\salto
\vskip2pt

   {\bf 6.2 Special features in  $ \Bfuno $  and  $ H^{(0)}_{f,f/2}
\; $.} \  The choice of parameter  $ \, x = 1 \, $  has two important
consequences.  Namely, the set  $ D_{f,f/2} $  of  $ \big(f/2\big) $--arc
$ f $--diagrams  is just a (multiplicative)  {\sl submonoid\/}  of
$ \Bfuno $,  and the action of  $ \Bfuno $  onto  $ H^{(0)}_{f,f/2} $
restricts to an action of  $ D_{f,f/2} $  onto the set  $ J_{f,f/2} $
of  $ \big(f,f/2\big) $--junctions.  Both facts follow directly from
definitions, which in fact give even more precise results, as follows.

\vskip4pt

   First observe that every  $ \, \d \in D_{f,f/2} \, $  is uniquely
determined by its arc structure  $ \as(\d) \, $  (notation of \S 1.2):
in short we write  $ \, \d \cong \as(\d) \, $.  Second, let  $ \, \d_1
\, $,  $ \, \d_2 \in D_{f,f/2} \, $.  Then by \S 1.3 one has  $ \, \d_1
\d_2 = \d_1 \ast \d_2 \in D_{f,f/2} \; \Big(\! \subset \Bfuno\big(f/2\big)
= \Bfuno\big[f/2\big] \Big) \, $,  \, and this product is uniquely
characterized by  $ \, \d_1 \d_2 \cong \as(\d_1 \d_2) =
\big( \tas(\d_1), \bas(\d_2) \big) \, $.

\vskip4pt

   Third, let  $ \, \d \in D_{f,f/2} \; \Big(\! \subset \Bfuno\big(f/2\big)
\Big) \, $  and  $ \, j \in J_{f,f/2} \; \Big(\! \subset H^{(0)}_{f,f/2}
\Big) \, $.  Then by definition (\S 2.6) the  $ \Bfuno $--action  on
$ H^{(0)}_{f,f/2} $  yields  $ \; \d.j = \tas(\d) \in J_{f,f/2} \;
\Big(\! \subset H^{(0)}_{f,f/2} \Big) \; $.

\vskip4pt

   We still need two more tools: the unique  $ \Bbbk $--linear  map
$ \; \text{\it Tr}_{\scriptscriptstyle \Cal{B}} \, \colon \, \Bfuno
\big(f/2\big) \longrightarrow \Bbbk \; $  such that  $ \, \text{\it
Tr}_{\scriptscriptstyle \Cal{B}}(\d) = 1 \, $  for all  $ \; \d \in
D_{f,f/2} \; $,  \, and the unique  $ \Bbbk $--linear  map  $ \; \text{\it
Tr}_{\scriptscriptstyle H} \, \colon \, H^{(0)}_{f,f/2} \longrightarrow
\Bbbk \; $  such that  $ \, \text{\it Tr}_{\scriptscriptstyle H}(j)
= 1 \, $  for all  $ \; j \in J_{f,f/2} \; $  (where notation  {\it Tr\/}
should be a reminder for ``trace'').

\salto
 \vskip2pt

\proclaim{Theorem 6.3}  Let  $ \; R^{(1)}_f := \, \text{\sl
$ \Bbbk $--span}\Big( \text{\sl Min}^{\,(1)}_{f;2} \, \bigcap \,
\Bfuno \! \big(f/2\big) \Big) \; $  like in  \S 5.8{\it (a)}.
Then
  $$  \Rad \left( \Bfuno\big(f/2\big) \right) \; = \; R^{(1)}_f \;
= \; \Ker\,\big(\text{\it Tr}_{\scriptscriptstyle B}\big) \; = \;
\text{$ \Bbbk $--{\sl span}}\Big(\big\{\, \d - \d' \,\big|\, \d,
\d' \in D_{f,f/2} \,\big\}\Big)   \leqno \indent \text{(a)}  $$
hence in particular Conjecture 5.10 holds true in this case  ($ \,
n = 1 \, $,  $ \, f \in 2 \N_+ \, $).
                                                       \par
   In addition, the semisimple quotient  $ \, {\Cal S}_f^{(1)}\big[f/2\big] $
of  $ \, \Bfuno\big[f/2\big] $  is simple of dimension 1.
  $$  \Rad\,\Big(H^{(0)}_{f,f/2}\Big) \,\; = \; R^{(1)}_f.\,H^{(0)}_{f,f/2}
\; = \; \Ker\,\big(\text{\it Tr}_{\scriptscriptstyle H}\big) \; = \;
\text{$ \Bbbk $--{\sl span}}\Big(\big\{\, j - j' \,\big|\, j, j' \in
J_{f,f/2} \,\big\}\Big)   \leqno \indent \text{(b)}  $$
In particular, Conjecture 5.14 is true for  $ \, n = 1 \, $,  $ \, f \in
2 \N_+ \, $  and  $ \, (k \, , \, \mu) = \big( f/2 \, , \, (0) \big) \, $.
                                                       \par
   In addition, the semisimple quotient of  $ \, H^{(0)}_{f,f/2} $  is
simple of dimension 1.

\endproclaim

\demo{Proof} {\it (a)} \,  Let us write  $ \, N(f) := \big|D_{f,f/2}\big|
= 2 (f-1)!! \, $,  and let  $ \, \d_1 , \, \d_2 \, , \, \dots \, , \,
\d_{N(f)} \, $  be a numbering of the elements of  $ \, D_{f,f/2} \; $.
It is clear that  $ \, \big\{\, \d_i - \d_{i+1} \, \big| \, i = 1, 2,
\dots, N(f) \! - \! 1 \big\} \, $  is a  $ \Bbbk $--basis  of  $ \, \Ker \, \big(
\text{\it Tr}_{\scriptscriptstyle \Cal{B}}\big) \, $,  \, and also that
$ \, \Ker \, \big( \text{\it Tr}_{\scriptscriptstyle \Cal{B}} \big) =
\text{$ \Bbbk $--{\sl span}} \Big( \big\{\, \d - \d' \,\big|\, \d,
\d' \in D_{f,f/2} \,\big\} \Big) \, $.
                                                      \par
   Now let  $ \, \d , \d' \in D_{f,f/2} \, $.  Then there exist
permutations  $ \, \sigma_+ , \sigma_- \in S_f \, $  such that  $ \;
\tas\big(\d'\big) = \d_{\sigma_+}.\,\tas(\d) \; $  and  $ \; \bas \big(
\d' \big) = \d_{\sigma_-}.\,\bas(\d) \; $.  Let us set  $ \, \d^* :=
\d_{\sigma_+} \! \cdot \d \, $   --- so that  $ \, \d^* \cong \big(
\d_{\sigma_+}.\,\tas(\d) \, , \bas(\d) \big) \, $  ---   let  $ \;
\sigma_\pm = \big( h_1^\pm \, k_1^\pm \big) \, \big( h_2^\pm \, k_2^\pm
\big) \cdots \big( h_{\ell(\sigma^\pm)} \, k_{\ell(\sigma^\pm)} \big)
\; $  be a factorisation of  $ \sigma_\pm $  into a product of
transpositions, and define
  $$  \displaylines{
   \d_0 := \d \; ,  \quad  \d_s := \d_{(h_{\ell(\sigma_+)-s+1}\,
k_{\ell(\sigma_+)-s+1})}.\,\d_{s-1}  \qquad  \big(\, s = 1, 2,
\dots, \ell(\sigma_+) \,\big)  \cr
   \d^*_0 := \d^* \; ,  \quad  \d^*_s := \d_{(h_{\ell(\sigma_-)-s+1}\,
k_{\ell(\sigma_-)-s+1})}.\,\d^*_{s-1}  \qquad  \big(\, s = 1, 2, \dots,
\ell(\sigma_-) \,\big)  \cr }  $$
(in particular,  $ \, \d_{\ell(\sigma_+)} = \d_{\sigma_+}.\,\d =: \d^*
=: \d^*_0 \; $).  Then we can expand  $ \, \d - \d' \, $  as
  $$  \d - \d'  \; = \;  {\textstyle \sum_{s=0}^{\ell(\sigma_+)-1}}
\big( \d_s - \d_{s+1} \big)  \, + \,  {\textstyle \sum_{s=0}^{\ell
(\sigma_-)-1}} \big( \d^*_s - \d^*_{s+1} \big)   \eqno (6.1)  $$
   \indent   By definition, any minor in  $ R^{(1)}_f $  is of the form
$ \, \d - \d' \, $,  with  $ \, \d, \d' \in D_{f,f/2} \, $  which differ
from each other only for a transposition of two bottom or two top vertices.
In other words, there are indices  $ \, h $,  $ k \in \{1,\dots,f\} \, $
such that   --- in notation of \S 1.3 ---
 we have\break
 \eject
\noindent
 $ \; \d' = \d \, \d_{(h\,k)}
\; $  (for  {\sl bottom\/}  vertices) or  $ \; \d' = \d_{(h\,k)} \, \d \; $
(for  {\sl top\/}  vertices).  Therefore, each summand in the right-hand
side of (6.1) above belongs to  $ \, R^{(1)}_f \, $,  \, hence  $ \; \big(
\d - \d' \big) \in R^{(1)}_f \; $  too.  By the previous analysis, it
follows that
 \vskip-7pt
  $$  R^{(1)}_f  \; = \;\,
\text{$ \Bbbk $--{\sl span}}\Big(\big\{\, \d - \d'
\,\big|\, \d, \d' \in D_{f,f/2} \,\big\}\Big)  \,\; = \;
\Ker\,\big(\text{\it Tr}_{\scriptscriptstyle \Cal{B}}\big)
\quad .  $$
 \vskip-2pt
   \indent   For every  $ \, \bar{\d} \, , \d \, , \d' \in D_{f,f/2} \, $,
\, we have also  $ \; \bar{\d} \, \big( \d - \d' \,\big) = \bar{\d} \, \d -
\bar{\d} \, \d' \in \Ker \, \big( \text{\it Tr}_{\scriptscriptstyle \Cal{B}}
\big) \, $,  \; so that  $ \Ker \, \big( \text{\it Tr}_{\scriptscriptstyle
\Cal{B}} \big) $  is a  $ \Bfuno\big(f/2\big) $--submodule  of  $ \, \Bfuno
\big(f/2\big) \, $  itself.  Then the quotient  $ \; \Bfuno\!\big(f/2\big)
\! \Big/ \! \Ker \, \big( \text{\it Tr}_{\scriptscriptstyle \Cal{B}} \big)
\, $,  which has dimension 1, is a  {\sl simple\/}  module for  $ \,
\Bfuno\!\big(f/2\big) = \Bfuno\!\big[f/2\big] $.
                                                      \par
   In addition, it is easy to see by direct computation that  $ \, \big(
\d_1 - \d'_1 \big) \big( \d_2 - \d'_2 \big) \big( \d_3 - \d'_3 \big) = 0
\, $  for every  $ \, \d_i \, , \d'_i \in D_{f,f/2} \; (i=1,2,3) \, $.
This implies that
  $$  \Big( \Ker \, \big( \text{\it Tr}_{\scriptscriptstyle \Cal{B}} \big)
\! \Big)^3  \, = \; \Big( R^{(1)}_f \Big)^3  \, = \;
  \Big( \text{$ \Bbbk $--{\sl span}} \Big( \big\{\,
\d - \d' \,\big|\, \d, \d' \in D_{f,f/2} \,\big\} \Big) \Big)^3  \, =
\;\, 0  $$
so that  $ \, R^{(1)}_f \! = \Ker \, \big( \text{\it
Tr}_{\scriptscriptstyle \Cal{B}} \big) \, $  is contained
in  $ \, \Rad\,\Big(\Bfuno\!\big(f/2\big)\Big) \, $.  But then
$ \; \Bfuno\!\big(f/2\big) \Big/ \! \Ker \, \big( \text{\it
Tr}_{\scriptscriptstyle \Cal{B}} \big) \, $  must be isomorphic to the
semisimple quotient  $ \, {\Cal S}_f^{(1)}\!\big[f/2\big] \, $  of  $ \,
\Bfuno\!\big[f/2\big] = \Bfuno\!\big(f/2\big) \, $;  \, therefore  $ \,
{\Cal S}_f^{(1)} \! \big[ f/2 \big] $  itself is simple of dimension 1,
and
  $$  \Rad \left( \Bfuno\big(f/2\big) \right)  \; = \;
\Ker\,\big(\text{\it Tr}_{\scriptscriptstyle B}\big)  \,\; = \; \,
\text{$ \Bbbk $--{\sl span}}\Big(\big\{\, \d - \d' \,\big|\, \d,
\d' \in D_{f,f/2} \,\big\}\Big)  \; = \;  R^{(1)}_f \; ,  \quad
\text{q.e.d.}  $$
   \indent   {\it (b)} \,  Let us write  $ \, J_{f,f/2} = \big\{\, j_1 , \,
j_2 \, , \, \dots \, , \, j_{n(f)} \,\big\} \, $,  \, with  $ \, n(f) :=
\big| J_{f,f/2} \big| = (f-1)!! \, $.  Then, like in the case of  $ \,
\Ker\,\big(\text{\it Tr}_{\scriptscriptstyle \Cal{B}}\big) \, $,  \, it
is clear that  $ \, \big\{\, j_s - j_{s+1} \,\big|\, s = 1, 2, \dots,
n(f)-1 \big\} \, $  is a  $ \Bbbk $--basis  of  $ \, \Ker \, \big(
\text{\it Tr}_{\scriptscriptstyle H} \big) \, $,  \, and  $ \, \Ker \,
\big( \text{\it Tr}_{\scriptscriptstyle H}\big) = \text{$ \Bbbk $--{\sl
span}}\Big(\big\{\, j - j' \,\big|\, j, j' \in J_{f,f/2} \,\big\}\Big)
\, $.  Since  $ \, \d.j = \tas(\d) \, $  for all  $ \, \d \in D_{f,f/2}
\, $  and all  $ \, j \in J_{f,f/2} \, $  (see \S 6.2), we have also
$ \, \d.\big( j - j' \big) = 0 \, $  (for all  $ \; j, j' \in J_{f,f/2}
\, $)  and so  $ \, \Ker \, \big( \text{\it Tr}_{\scriptscriptstyle H}
\big) \, $  is a  $ \Bfuno\!\big(f/2\big) $--submodule  of  $ \, H^{(0)}_{f,
f/2} \; $.  As the quotient  $ \; H^{(0)}_{f,f/2} \Big/ \! \Ker \, \big(
\text{\it Tr}_{\scriptscriptstyle H} \big) \; $  is  $ 1 $--dimensional,
we conclude that it is a  {\sl simple\/}  $ \Bfuno\!\big(f/2\big) $--module.
                                                      \par
   Now pick any  $ \, j $,  $ j' \in J_{f,f/2} \, $:  \, then there exists
a permutation  $ \, \sigma \in S_f \, $  such that  $ \, j' = \d_\sigma.j
\, $  (notation of \S 1.3).  Let  $ \; \sigma = (h_1\,k_1) \, (h_2\,k_2)
\cdots (h_{\ell(\sigma)}\,k_{\ell(\sigma)}) \; $  be a factorisation of
$ \sigma $  into a product of transpositions, and let  $ \, j_0 := j \, $,
$ \, j_s := \d_{(h_{\ell(\sigma)-s+1} \, k_{\ell(\sigma)-s+1})}.j_{s-1}
\, $  for all  $ \, s = 1, 2, \dots, \ell(\sigma) \, $.  Then we can
expand  $ \, j - j' \, $  as  $ \; j - j' \, = \, j - \d_\sigma.j \,
= \, \sum_{s=0}^{\ell(\sigma)-1} \! \big( j_s - j_{s+1} \big) \; $.
                                                      \par
   Now take any  $ \, y \in J_{f,f/2} \, $  and let  $ \, \d_s \in
D_{f,f/2} \, $  be characterized by  $ \, \as(\d_s) \cong \big( j_s
\, , y \big) \, $.  Recalling that  $ \, \d.y' = \tas(\d) \, $,  for
all  $ \, \d \in D_{f,f/2} \, $  and  $ \, y' \in J_{f,f/2} \, $
(see \S 6.2), we have
  $$  \Big({\textstyle \sum_{s=0}^{\ell(\sigma)-1}} \big( \d_s - \d_{s+1}
\big)\Big).\,y \, = \, {\textstyle \sum_{s=0}^{\ell(\sigma)-1}}
\big( \tas(\d_s) - \tas(\d_{s+1}) \big) \, = \, {\textstyle
\sum_{s=0}^{\ell(\sigma)-1}} \big( j_s - j_{s+1} \big) \,
= \, j - j' \, .  $$
By the previous description of  $ R^{(1)}_f $  (see part  {\it (a)\/})  we
have  $ \, {\textstyle \sum_{s=0}^{\ell(\sigma)-1}} \big( \d_s - \d_{s+1}
\big) \in R^{(1)}_f \, $,  \, hence we can conclude that  $ \; \big( j -
j' \big) \in R^{(1)}_f.\,H^{(0)}_{f,f/2} \; $.  As we have already shown
that  $ \; \text{$ \Bbbk $--{\sl span}}\Big(\big\{\, j - j' \,\big|\,
j, j' \in J_{f,f/2} \,\big\}\Big) = \Ker \, \big( \text{\it
Tr}_{\scriptscriptstyle H} \big) \; $  and  $ \; R^{(1)}_f =
\Rad\,\Big(\Bfuno\!\big(f/2\big)\!\Big) \, $,  \, we have
 \vskip-3pt
  $$  \Ker \, \big( \text{\it Tr}_{\scriptscriptstyle H} \big)  \,\; = \;
\Rad \, \Big( \Bfuno \! \big(f/2\big) \! \Big).\,H^{(0)}_{f,f/2}  \,\;
\subseteq \;  \Rad\,\Big(H^{(0)}_{f,f/2}\Big)  $$
 \vskip-1pt
\noindent
 where the last inclusion follows from the standard inclusion  $ \, \Rad\,
(A).M \subseteq \Rad\,(M) \, $,  \; which holds for every ring  $ A $
and every  $ A $--module  $ M $.
 Finally, as we saw that the quo-\break
 \eject
\noindent
 tient  $ \; H^{(0)}_{f,f/2} \Big/ \! \Ker \, \big( \text{\it
Tr}_{\scriptscriptstyle H} \big) \; $  is a  ($ 1 $--dimensional)
{\sl simple\/}  $ \Bfuno\!\big(f/2\big) $--module,  we can conclude
that  $ \; \Ker \, \big( \text{\it Tr}_{\scriptscriptstyle H} \big)
\supseteq \Rad \, \Big( H^{(0)}_{f,f/2} \Big) \; $  as well, thus
eventually  $ \; \Ker \, \big( \text{\it Tr}_{\scriptscriptstyle H}
\big) = R^{(1)}_f.\,H^{(0)}_{f,f/2} = \Rad\,\Big(H^{(0)}_{f,f/2}\Big)
\, $,  \; and the semisimple quotient of  $ H^{(0)}_{f,f/2} $  is
simple of dimension 1, q.e.d.   \qed
\enddemo

\salto
% \vskip1pt

\vskip0,31truecm

{}

\enddocument